\documentclass[12pt]{amsart}
\usepackage{amssymb,amsthm,amsmath, ulem, fullpage}
\usepackage{alltt,calc}
\usepackage[all]{xy}

\DeclareFontFamily{OMS}{rsfs}{\skewchar\font'60}
\DeclareFontShape{OMS}{rsfs}{m}{n}{<-5>rsfs5 <5-7>rsfs7 <7->rsfs10 }{}
\DeclareSymbolFont{rsfs}{OMS}{rsfs}{m}{n}
\DeclareSymbolFontAlphabet{\scr}{rsfs}
\newtheorem{theorem}{Theorem}[section]
\newtheorem{lemma}[theorem]{Lemma}
\newtheorem{proposition}[theorem]{Proposition}
\newtheorem{corollary}[theorem]{Corollary}
\newtheorem{conjecture}[theorem]{Conjecture}
\theoremstyle{definition}
\newtheorem{definition}[theorem]{Definition}
\newtheorem{example}[theorem]{Example}

\theoremstyle{remark}
\newtheorem{remark}[theorem]{Remark}

\newtheorem{question}[theorem]{Question}

\newcommand{\myR}{{R}}
\newcommand{\myH}{{\bf \mathcal H}}
\newcommand{\tld}{\widetilde }
\newcommand{\mydot}{\centerdot}
\newcommand{\blank}{\underline{\hskip 10pt}}
\newcommand{\DuBois}[1]{{\uline \Omega {}^0_{#1}}}
\newcommand{\FullDuBois}[1]{{\uline \Omega {}^{\mydot}_{#1}}}

\newcommand{\uSch}{\underline{\kSch}}

\newcommand{\qis}{\simeq_{\text{qis}}}

\newcommand{\sF}{\scr{F}}

\newcommand{\sH}{\scr{H}}
\newcommand{\sI}{\scr{I}}

\newcommand{\sL}{\scr{L}}

\newcommand{\sO}{\scr{O}}

\newcommand{\mJ}{\mathcal{J}}
\newcommand{\uTwo}{\underline{2}}
\newcommand{\uOne}{\underline{1}}

\newcommand{\bm}{\mathfrak{m}}
\newcommand{\bA}{\mathbb{A}}

\newcommand{\bC}{\mathbb{C}}

\newcommand{\bH}{\mathbb{H}}

\newcommand{\bN}{\mathbb{N}}

\newcommand{\bP}{\mathbb{P}}
\newcommand{\bQ}{\mathbb{Q}}

\newcommand{\bZ}{\mathbb{Z}}
\newcommand{\al}{\alpha}
\newcommand{\be}{\beta}
\newcommand{\ga}{\gamma}

\renewcommand{\O}{\scr O}

\newcommand{\N}{\mathbb {N}}

\renewcommand{\O}{\mbox{$\mathcal{O}$}}

\newcommand{\tensor}{\otimes}

\DeclareMathOperator{\rd}{{rd}}
\DeclareMathOperator{\an}{{an}}

\DeclareMathOperator{\coherent}{{coh}}

\DeclareMathOperator{\sn}{{sn}}

\DeclareMathOperator{\id}{{id}}

\DeclareMathOperator{\opp}{{op}}

\DeclareMathOperator{\ob}{{Ob}}

\DeclareMathOperator{\Supp}{{Supp}}
\DeclareMathOperator{\Sing}{{Sing}}

\DeclareMathOperator{\Spec}{{Spec}}

\DeclareMathOperator{\Hom}{Hom}

\DeclareMathOperator{\red}{red}

\def\spec#1.#2.{{\bold S\bold p\bold e\bold c}_{#1}#2}
\def\proj#1.#2.{{\bold P\bold r\bold o\bold j}_{#1}\sum #2}
\def\ring#1.{\scr O_{#1}}
\def\map#1.#2.{#1 \to #2}
\def\longmap#1.#2.{#1 \longrightarrow #2}
\def\factor#1.#2.{\left. \raise 2pt\hbox{$#1$} \right/
\hskip -2pt\raise -2pt\hbox{$#2$}}
\def\pe#1.{\mathbb P(#1)}
\def\pr#1.{\mathbb P^{#1}}

\def\coh#1.#2.#3.{H^{#1}(#2,#3)}
\def\dimcoh#1.#2.#3.{h^{#1}(#2,#3)}
\def\hypcoh#1.#2.#3.{\mathbb H_{\vphantom{l}}^{#1}(#2,#3)}
\def\loccoh#1.#2.#3.#4.{H^{#1}_{#2}(#3,#4)}
\def\dimloccoh#1.#2.#3.#4.{h^{#1}_{#2}(#3,#4)}
\def\lochypcoh#1.#2.#3.#4.{\mathbb H^{#1}_{#2}(#3,#4)}
\def\ses#1.#2.#3.{0  \longrightarrow  #1   \longrightarrow
 #2 \longrightarrow #3 \longrightarrow 0}
\def\sesshort#1.#2.#3.{0
 \rightarrow #1 \rightarrow #2 \rightarrow #3 \rightarrow 0}
 
\def\iff#1#2#3{
    \hfil\hbox{\hsize =#1 \vtop{\noin #2} \hskip.5cm
    \lower.5\baselineskip\hbox{$\Leftrightarrow$}\hskip.5cm
    \vtop{\noin #3}}\hfil\medskip}
\def\myoplus#1.#2.{\underset #1 \to {\overset #2 \to \oplus}}

\newcounter{are-there-sections}
\makeatletter

\renewcommand\subsection{
  \renewcommand{\sfdefault}{pag}
  \@startsection{subsection}%
  {2}{0pt}{-\baselineskip}{.2\baselineskip}{\raggedright
    \sffamily\itshape\small
  }}
\renewcommand\section{
  \renewcommand{\sfdefault}{phv}
  \@startsection{section} %
  {1}{0pt}{\baselineskip}{.2\baselineskip}{\centering
    \sffamily
    \scshape
}}
\newcounter{lastyear}
\addtocounter{lastyear}{-1}

\newcommand\noin{\noindent}

\newtheoremstyle{bozont}{3pt}{3pt}%
     {\itshape}
     {}
     {\bfseries}
     {.}
     {.5em}
     {\thmname{#1}\thmnumber{ #2}\thmnote{ \rm #3}}
\newtheoremstyle{bozont-sf}{3pt}{3pt}%
     {\itshape}
     {}
     {\sffamily}
     {.}
     {.5em}
     {\thmname{#1}\thmnumber{ #2}\thmnote{ \rm #3}}
\newtheoremstyle{bozont-sc}{3pt}{3pt}%
     {\itshape}
     {}
     {\scshape}
     {.}
     {.5em}
     {\thmname{#1}\thmnumber{ #2}\thmnote{ \rm #3}}
\newtheoremstyle{bozont-remark}{3pt}{3pt}%
     {}
     {}
     {\scshape}
     {.}
     {.5em}
     {\thmname{#1}\thmnumber{ #2}\thmnote{ \rm #3}}
\newtheoremstyle{bozont-def}{3pt}{3pt}%
     {}
     {}
     {\bfseries}
     {.}
     {.5em}
     {\thmname{#1}\thmnumber{ #2}\thmnote{ \rm #3}}
\newtheoremstyle{bozont-reverse}{3pt}{3pt}%
     {\itshape}
     {}
     {\bfseries}
     {.}
     {.5em}
     {\thmnumber{#2.}\thmname{ #1}\thmnote{ \rm #3}}
\newtheoremstyle{bozont-reverse-sc}{3pt}{3pt}%
     {\itshape}
     {}
     {\scshape}
     {.}
     {.5em}
     {\thmnumber{#2.}\thmname{ #1}\thmnote{ \rm #3}}
\newtheoremstyle{bozont-reverse-sf}{3pt}{3pt}%
     {\itshape}
     {}
     {\sffamily}
     {.}
     {.5em}
     {\thmnumber{#2.}\thmname{ #1}\thmnote{ \rm #3}}
\newtheoremstyle{bozont-remark-reverse}{3pt}{3pt}%
     {}
     {}
     {\sc}
     {.}
     {.5em}
     {\thmnumber{#2.}\thmname{ #1}\thmnote{ \rm #3}}
\newtheoremstyle{bozont-def-reverse}{3pt}{3pt}%
     {}
     {}
     {\bfseries}
     {.}
     {.5em}
     {\thmnumber{#2.}\thmname{ #1}\thmnote{ \rm #3}}
\newtheoremstyle{bozont-def-newnum-reverse}{3pt}{3pt}%
     {}
     {}
     {\bfseries}
     {}
     {.5em}
     {\thmnumber{#2.}\thmname{ #1}\thmnote{ \rm #3}}
\theoremstyle{bozont}
\ifnum \value{are-there-sections}=0 {%
  \newtheorem{proclaim}{Theorem}
}
\else {%
  
}
\fi
\newtheorem{thm}[theorem]{Theorem}

\newtheorem{conj}[theorem]{Conjecture}

\theoremstyle{bozont-sc}
\newtheorem{theorem-special}[theorem]{\specialthmname}

\theoremstyle{bozont-remark}

\newtheorem{subrem}[equation]{Remark}

\newtheorem*{SubHeading*}{\SubHeadingName}%
\newtheorem{SubHeading}[theorem]{\SubHeadingName}
\newtheorem{sSubHeading}[equation]{\sSubHeadingName}
\newenvironment{demo}[1] {\def\SubHeadingName{#1}\begin{SubHeading}}
  {\end{SubHeading}}%
\newenvironment{demo-r}[1]{\def\SubHeadingName{#1}\begin{SubHeading-r}}
  {\end{SubHeading-r}}%
\newenvironment{subdemo-r}[1]{\def\sSubHeadingName{#1}\begin{sSubHeading-r}}
  {\end{sSubHeading-r}} %
\newenvironment{demo*}[1]{\def\SubHeadingName{#1}\begin{SubHeading*}}
  {\end{SubHeading*}}%
\newtheorem{defini}[theorem]{Definition}

\newtheorem*{ack}{Acknowledgements}

\newtheorem{defn-thm}[theorem]{Definition--Theorem}  
\theoremstyle{bozont-def}
\newtheorem{defn}[theorem]{Definition}

\theoremstyle{bozont-reverse}

\theoremstyle{bozont-reverse-sc}
\newtheorem{theoremr-special}[theorem]{\specialthmname}
{\def\specialthmname{#1}\begin{theoremr-special}}%
{\end{theoremr-special}}
\theoremstyle{bozont-remark-reverse}

\newtheorem{SubHeading-r}[theorem]{\SubHeadingName}
\newtheorem{sSubHeading-r}[equation]{\sSubHeadingName}
\newtheorem{SubHeadingr}[theorem]{\SubHeadingName}

\theoremstyle{bozont-def-newnum-reverse}

\theoremstyle{bozont-def-reverse}

\newtheorem{newnumspecial}[theorem]{\specialnewnumname}

\numberwithin{equation}{theorem}
\numberwithin{figure}{section}

\newcounter{firstsubsection} %
\newcommand\Subsection[1]
      { %
        \ifnum \value{firstsubsection}=0 %
        \setcounter{subsection}{\value{theorem}} %
        \setcounter{firstsubsection}{1} %
        \fi %
        \noin\subsection{#1} %
        \numberwithin{theorem}{subsection}
        \numberwithin{equation}{subsection}
        \numberwithin{figure}{subsection} %
        }

\newcounter{rosternumber}

\newenvironment{enumerate-p}{
  \begin{enumerate}}
  {\setcounter{equation}{\value{enumi}}\end{enumerate}}
\newenvironment{enumerate-cont}{
  \begin{enumerate}
    {\setcounter{enumi}{\value{equation}}}}
  {\setcounter{equation}{\value{enumi}}
  \end{enumerate}}

\newlength{\swidth}
\makeatother

\newcounter{append}%
\renewcommand\appendix[1]{%
  \addtocounter{append}{1}
  \renewcommand\thesection{Appendix \Alph{append}}%
  \section{#1}
  \renewcommand\thesection{\Alph{append}}%
}

\DeclareMathAlphabet{\smallchanc}{OT1}{pzc}%
                                 {m}{it}
\DeclareFontFamily{OT1}{pzc}{}
\DeclareFontShape{OT1}{pzc}{m}{it}%
             {<-> s * [1.100] pzcmi7t}{}
\DeclareMathAlphabet{\mathchanc}{OT1}{pzc}%
                                 {m}{it}

\newcommand{\mcH}{\mathchanc{H}}

\newcommand{\mcm}{\mathchanc{m}}

\newcommand{\mco}{\mathchanc{o}}

\DeclareFontFamily{OMS}{rsfs}{\skewchar\font'60}
\DeclareFontShape{OMS}{rsfs}{m}{n}{<-5>rsfs5 <5-7>rsfs7 <7->rsfs10 }{}
\DeclareSymbolFont{rsfs}{OMS}{rsfs}{m}{n}
\DeclareSymbolFontAlphabet{\scr}{rsfs}

\newcommand{\sfC}{{\sf C}}
\newcommand{\sfD}{{\sf D}}

\newcommand{\ff}{\mathfrak{f}}

\newcommand{\into}{\hookrightarrow}

\newcommand{\wtilde}{\widetilde}
\newcommand{\wt}{\widetilde}

\newcommand{\rdown}[1]{\lfloor{#1}\rfloor}

\newcommand{\leteq}{\colon\!\!\!=}
\newcommand{\col}{\colon}

\DeclareMathOperator{\opdiv}{div}
\DeclareMathOperator{\Exc}{Exc}
\newcommand{\sHom}[0]{{\mcH\mco\mcm}}

\DeclareMathOperator{\Ob}{Ob}

\DeclareMathOperator{\reg}{reg}

\newcommand{\cx}[1]{{#1}^{\raisebox{.15em}{\ensuremath\centerdot}}}
\newcommand{\Om}{\underline{\Omega}}

\def\coh#1.#2.#3.{H^{#1}(#2,#3)}
\def\dimcoh#1.#2.#3.{h^{#1}(#2,#3)}
\def\hypcoh#1.#2.#3.{\mathbb H_{\vphantom{l}}^{#1}(#2,#3)}
\def\loccoh#1.#2.#3.#4.{H^{#1}_{#2}(#3,#4)}
\def\dimloccoh#1.#2.#3.#4.{h^{#1}_{#2}(#3,#4)}
\def\lochypcoh#1.#2.#3.#4.{\mathbb H^{#1}_{#2}(#3,#4)}
\def\dist#1.#2.#3.{  #1   \longrightarrow
 #2 \longrightarrow #3 \stackrel{+1}{\longrightarrow} } 
\def\CDdist#1.#2.#3.{  #1   @>>>  #2  @>>>   #3 @>+1>> }
\def\shortses#1.#2.#3.{0  \rightarrow  #1   \rightarrow
 #2  \rightarrow   #3 \rightarrow  0}
\def\shortdist#1.#2.#3.{  #1   \rightarrow
 #2  \rightarrow   #3 \stackrel{+1}{\rightarrow} }  
\def\ddist#1.#2.#3.#4.#5.#6.{\CD
#1 @>>> #2 @>>> #3 @>+1>> \\
@VVV @VVV @VVV \\
#4 @>>> #5 @>>> #6 @>+1>>
\endCD}
\def\ddistun#1.#2.#3.#4.#5.#6.{\CD
#1 @>>> #2 @>>> #3 @>+1>> \\
@. @VVV @VVV  \\
#4 @>>> #5 @>>> #6 @>+1>>
\endCD}
\def\Iff#1#2#3{
\hfil\hbox{\hsize =#1
\vtop{\noin #2}
\hskip.5cm
\lower.5\baselineskip\hbox{$\Leftrightarrow$}\hskip.5cm
\vtop{\noin #3}}\hfil\medskip}
\newcommand{\union}\cup
\newcommand{\intersect}\cap
\newcommand{\Union}\bigcup
\newcommand{\Intersect}\bigcap
\newcommand{\resto}{\big\vert_}

\def\qis{\,{\simeq}_{qis}\,}

\def\ww#1.#2.{\curlywedge
_{#1}^{#2}}
\def\wws#1.#2.{\scriptstyle\ww#1.#2.}

\def\fh#1.#2.{F^{#1} R^{#2}}
\def\epq#1.#2.#3.{E^{#1,#2}_{#3}}
\def\dd#1.#2.#3.{d^{\, #1,#2}_{#3}}
\def\kk#1.#2.#3.{K^{\, #1,#2}_{#3}}
\def\ki#1.#2.#3.{I^{\, #1,#2}_{#3}}
\def\ff#1.#2.{\mathfrak F_{#1}^{#2}}
\def\ffb#1.#2.{\boxed{\mathfrak F_{#1}^{#2}}}
\def\ffbx#1.#2.#3.{\boxed{\mathfrak F_{#1}^{#2}(#3)}}
\def\ffs#1.#2.{\scriptstyle\mathfrak F_{#1}^{#2}}
\def\ffsb#1.#2.{\boxed{\ffs#1.#2.}}
\def\ffsbx#1.#2.#3.{\boxed{\scriptstyle\mathfrak F_{#1}^{#2}(#3)}}
\def\fa#1.#2.{\mathfrak f_{#1}^{#2}}
\def\fab#1.#2.{\boxed{\mathfrak f_{#1}^{#2}}}
\def\ffa#1.#2.{\mathfrak F\mathfrak i\mathfrak l\mathfrak t_{#1}^{#2}}
\def\ffab#1.#2.{\boxed{\mathfrak F\mathfrak i\mathfrak l\mathfrak t_{#1}^{#2}}}
\def\rfi#1.#2.{R^{#1}\functor(#2\,)}

\def\al#1.#2.#3.{\alpha_{#1,#2}^{#3}}
\def\be#1.#2.{\beta_{#1,#2}}
\def\ga#1.#2.{\gamma_{#1,#2}}
\def\yy#1.#2.#3.#4.{y_{#1,#2}^{#3}(#4)}
\def\zz#1.#2.#3.#4.{z_{#1,#2}^{#3}(#4)}

\renewcommand{\O}{\sO}
\renewcommand{\uSch}{\sf Sch_{\text{red}}}
\DeclareMathOperator{\Mor}{{Mor}}

\newcommand{\uFaisc}{{\sf Sheaves}}
\renewcommand{\myH}{h}
\renewcommand{\mydot}{{{\,\begin{picture}(1,1)(-1,-2)\circle*{2}\end{picture}\ }}}
\numberwithin{equation}{theorem}
\begin{document}

\title{Hodge theory meets the minimal model program: a survey of log canonical and
  Du~Bois singularities}

\author{S\'andor J Kov\'acs and Karl Schwede}

\address{\vskip-.8cm\ \newline \noindent S\'andor J Kov\'acs: \sf University of
  Washington, Department of Mathematics, Seattle, WA 98195, USA}

\email{kovacs@math.washington.edu}

\address{\vskip -.65cm \noindent Karl E.\ Schwede: \sf Department of Mathematics,
  University of Michigan.  Ann Arbor, Michigan 48109-1109}

\email{kschwede@umich.edu}

\date{\today}


\thanks{The first named author was supported in part by NSF Grant DMS-0554697 and the
  Craig McKibben and Sarah Merner Endowed Professorship in Mathematics.}

\thanks{The second named author was partially supported by RTG grant number 0502170
  and by a National Science Foundation Postdoctoral Research Fellowship.}

\begin{abstract}
  This is a survey of some recent developments in the study of singularities related
  to the classification theory of algebraic varieties. In particular, the definition
  and basic properties of {Du~Bois singularities} and their connections to the more
  commonly known singularities of the minimal model program are reviewed and
  discussed.
\end{abstract}

\maketitle

\section{Introduction}The primary goal of this note is to survey some recent
developments in the study of singularities related to the minimal model program. In
particular, we review the definition and basic properties of \emph{Du~Bois
  singularities} and explain how these singularities fit into the minimal model
program and moduli theory.

Since we can resolve singularities \cite{Hironaka64}, one might ask the question why
we care about them at all.  It turns out that in various situations we are forced to
work with singularities even if we are only interested in understanding smooth
objects.

One reason we are led to study singular varieties is provided by the minimal model
program \cite{KollarMori}. The main goal is classification of algebraic varieties and
the plan is to find reasonably simple representatives of all birational classes and
then classify these representatives. It turns out that the simplest objects in a
birational class tend to be singular. What this really means is that when choosing a
birational representative, we aim to have simple \emph{global} properties and this is
often achieved by a singular variety. Being singular means that there are points
where the \emph{local} structure is more complicated than on a smooth variety, but
that allows for the possibility of still having a somewhat simpler global structure
and along with it, good local properties at most points.

Another reason to study singularities is that to understand smooth objects we should
also understand how smooth objects may deform and degenerate.  This leads to the need
to construct and understand moduli spaces. And not just moduli for the smooth
objects.  Degenerations provide important information as well. In other words, it is
always useful to work with complete moduli problems, i.e., extend our moduli functor
so it admits a compact (and preferably projective) coarse moduli space.  This also
leads to having to consider singular varieties.

On the other hand, we have to be careful to limit the kind of singularities that we
allow in order to be able to handle them. One might view this survey as a list of the
singularities that we must deal with to achieve the above stated goals. Fortunately,
it is also a class of singularities with which we have a reasonable chance to be able
to work.

In particular, we will review Du~Bois singularities and related notions
including some very recent important results. We will also review a family of
singularities defined via characteristic $p$ methods, the Frobenius morphism, and
their connections to the other set of singularities we are discussing.

\newenvironment{nlist}%
{\begin{list}{}{\let\makelabel\nitem
      \setlength\labelwidth{0pt}
      \setlength\leftmargin{\labelwidth+\labelsep}}}
  {\end{list}
}
\newcommand*\nitem[1]{\textsf{#1$\bullet$}\hfil}

\begin{demo}  {Definitions and notation} %
  Let $k$ be an algebraically closed field.  
  Unless otherwise stated, all objects will be assumed to be defined over $k$. A
  \emph{scheme} will refer to a scheme of finite type over $k$ and unless stated
  otherwise, a \emph{point} refers to a closed point.

  For a morphism $Y\to S$ and another morphism $T\to S$, the symbol $Y_T$ will denote
  $Y\times_S T$. In particular, for $t\in S$ we write $X_t = f^{-1}(t)$. In addition,
  if $T=\Spec F$, then $Y_T$ will also be denoted by~$Y_F$.

  Let $X$ be a scheme and $\sF$ an $\sO_X$-module. The \emph{$m^\text{th}$ reflexive
    power} of $\sF$ is the double dual (or reflexive hull) of the $m^\text{th}$
  tensor power of $\sF$:
  $$
  \sF^{[m]}\leteq (\sF^{\otimes m})^{**}.
  $$
  A \emph{line bundle} on $X$ is an invertible $\sO_X$-module. A \emph{$\bQ$-line
    bundle} $\sL$ on $X$ is a reflexive $\sO_X$-module of rank $1$ that possesses a
  reflexive power which is a line bundle, i.e., there exists an $m\in \N_+$ such that
  $\sL^{[m]}$ is a line bundle.  The smallest such $m$ is called the \emph{index} of
  $\sL$.

  \begin{nlist}
  \item For the advanced reader: whenever we mention {Weil divisors}, assume that $X$
    is $S_2$ \cite[Thm.~8.22A(2)]{Hartshorne77} and think of a \emph{Weil divisorial
      sheaf}, that is, a rank $1$ reflexive $\sO_X$-module which is locally free in
    codimension $1$. For flatness issues consult \cite[Theorem
    2]{kollar-hulls-and-husks}.
  \item For the novice: whenever we mention Weil divisors, assume that $X$ is normal
    and adopt the definition \cite[p.130]{Hartshorne77}.
  \end{nlist}
  For a Weil divisor $D$ on $X$, its associated \emph{Weil divisorial sheaf} is the
  $\sO_X$-module $\sO_X(D)$ defined on the open set $U\subseteq X$ by the formula
  \begin{multline*}
    \Gamma (U, \sO_X(D))=\left\{ \frac ab \ \bigg|\ a,b\in \Gamma(U, \sO_X), b \text{
        is not a zero divisor anywhere on $U$, and } \right. \\ \left.
      \phantom{\bigg|} D|_U+\opdiv_U(a)-\opdiv_U(b)\geq 0 \right\}
  \end{multline*}
  and made into a sheaf by the natural restriction maps.

  A Weil divisor $D$ on $X$ is a \emph{Cartier divisor}, if its associated Weil
  divisorial sheaf, $\sO_X(D)$ is a line bundle. If the associated Weil divisorial
  sheaf, $\sO_X(D)$ is a $\bQ$-line bundle, then $D$ is a \emph{$\bQ$-Cartier
    divisor}. The latter is equivalent to the property that there exists an $m\in
  \bN_+$ such that $mD$ is a Cartier divisor.  Weil divisors form an abelian group.
  Tensoring this group with $\bQ$ (over $\bZ$) one obtains the group of
  \emph{$\bQ$-divisors} on $X$ (note that if $X$ is not normal, some unexpected
  things can happen in this process, see \cite[Chapter 16]{K+92}).

  The symbol $\sim$ stands for \emph{linear} and $\equiv$ for \emph{numerical
    equivalence} of divisors.

  Let $\sL$ be a line bundle on a scheme $X$. It is said to be \emph{generated by
    global sections} if for every point $x\in X$ there exists a global section
  $\sigma_x\in \coh 0.X.\sL.$ such that the germ $\sigma_x$ generates the stalk
  $\sL_x$ as an $\sO_X$-module. If $\sL$ is generated by global sections, then the
  global sections define a morphism
  $$
  \phi_{\sL}\col X\to\bP^N= \bP\left(\coh 0.X.\sL.^{*}\right).
  $$
  $\sL$ is called \emph{semi-ample} if $\sL^{m}$ is generated by global sections for
  $m\gg 0$. $\sL$ is called \emph{ample} if it is semi-ample and $\phi_{\sL^m}$ is an
  embedding for $m\gg 0$. A line bundle $\sL$ on $X$ is called \emph{big} if the
  global sections of $\sL^{m}$ define a rational map $\phi_{\sL^m}\col X\dasharrow
  \bP^N$ such that $X$ is birational to $\phi_{\sL^m}(X)$ for $m\gg 0$. Note that in
  this case $\sL^m$ need not be generated by global sections, so
  $\phi_{\sL^m}$ is not necessarily defined everywhere. We leave it for the reader
  the make the obvious adaptation of these notions for the case of $\bQ$-line
  bundles.

  The \emph{canonical divisor} of a scheme $X$ is denoted by $K_X$ and the
  \emph{canonical sheaf} of $X$ is denoted by $\omega_X$.

  A smooth projective variety $X$ is of \emph{general type} if $\omega_X$ is big.  It
  is easy to see that this condition is invariant under birational equivalence
  between smooth projective varieties. An arbitrary projective variety is of
  \emph{general type} if so is a desingularization of it.

  A projective variety is \emph{canonically polarized} if $\omega_X$ is ample. Notice
  that if a smooth projective variety is canonically polarized, then it is of general
  type.
\end{demo}

\begin{ack}
  We would like to thank Kevin Tucker, Zsolt Patakfalvi and the referee for reading a
  preliminary draft and making helpful suggestion for improving the presentation.
\end{ack}

\section{Pairs and resolutions}

For the reader's convenience, we recall a few definitions regarding pairs.

\begin{defn}\label{def:everythinglog1}
  A \emph{pair} $(X,\Delta)$ consists of a normal\footnote{Occasionally, we will
    discuss pairs in the non-normal setting.  See Section \ref{SubsectionSLC} for
    more details.} quasi-projective variety or complex space $X$ and an effective
  $\bQ$-divisor $\Delta \subset X$.  A \emph{morphism of pairs} $\gamma: (\wtilde X,
  \wtilde \Delta) \to (X,\Delta)$ is a morphism $\gamma: \wtilde X \to X$ such that
  $\gamma(\Supp(\wt\Delta) ) \subseteq \Supp(\Delta)$. A morphism of pairs $\gamma:
  (\wtilde X, \wtilde \Delta) \to (X,\Delta)$ is called \emph{birational} if it
  induces a birational morphism $\gamma:\wt X\overset\sim\to X$ and $\gamma(\wt
  \Delta)=\Delta$. It is an \emph{isomorphism} if it is birational and it induces an
  isomorphism $\gamma:\wt X\overset\simeq\to X$.
\end{defn}

\begin{defn}\label{def:everythinglog2a}
  Let $(X, \Delta)$ be a pair, and $x\in X$ a point. We say that $(X, \Delta)$ is
  \emph{snc at $x$}, if there exists a Zariski-open neighborhood $U$ of $x$ such that
  $U$ is smooth and $\Delta \cap U$ is reduced and has only simple normal crossings
  (see Section \ref{subsectionNormalCrossings} for additional discussion).  The pair
  $(X, \Delta)$ is \emph{snc} if it is snc at all $x\in X$.

  Given a pair $(X,\Delta)$, let $(X,\Delta)_{\reg}$ be the maximal open set of $X$
  where $(X,\Delta)$ is snc, and let $(X,\Delta)_{\Sing}$ be its complement, with the
  induced reduced subscheme structure.
\end{defn}

\begin{subrem} If a pair $(X, \Delta)$ is snc at a point $x$, this implies that all
  components of $\Delta$ are smooth at $x$. If instead of the condition that $U$ is
  Zariski-open one would only require this analytically locally, then
  Definition~\ref{def:everythinglog2a} would define normal crossing pairs rather than
  pairs with simple normal crossing.
\end{subrem}

\begin{defn}\label{def:everythinglog2b}
  A \emph{log resolution} of $(X, \Delta)$ is a proper birational morphism of pairs
  $\pi : (\tld X,\wt \Delta) \rightarrow (X,\Delta)$ that satisfies the following
  four conditions:
  \begin{enumerate-p}
  \item $\tld X$ is smooth,
  \item $\wt\Delta=\pi^{-1}_*\Delta$ is the strict transform of $\Delta$,
  \item $\Exc(\pi)$ is of pure codimension $1$,
  \item $\Supp(\tld \Delta \cup \Exc(\pi))$ is a simple normal crossings divisor.
  \end{enumerate-p}

  If in addition,
  \begin{enumerate-cont}
  \item the strict transform $\tld \Delta$ of $\Delta$ has smooth support,
  \end{enumerate-cont}
  then we call $\pi$ an \emph{embedded resolution} of $\Delta\subset X$.

  In many cases, it is also useful to require that $\pi$ is an isomorphism over $(X,
  \Delta)_{\reg}$.
\end{defn}






\section{Introduction to the singularities of the mmp}

Even though we have introduced pairs and most of these singularities make sense for
pairs, to make the introduction easier to digest we will mostly discuss the case when
$\Delta=\emptyset$.
As mentioned in the introduction, one of our goals is to show why we are forced to
work with singular varieties even if our primary interest lies with smooth varieties.

\subsection{Canonical singularities}

For an excellent introduction to this topic the reader is urged to take a thorough
look at Miles Reid's Young Person's Guide \cite{Reid87}. Here we will only touch on
the subject.

Let us suppose that we would like to get a handle on some varieties. Perhaps we want
to classify them or make some computations. In any case, a useful thing to do is to
embed the object in question into a projective space (if we can). Doing so requires a
(very) ample line bundle. It turns out that in practice these can be difficult to
find. In fact, it is not easy to find any non-trivial line bundle on an abstract
variety.

One possibility, when $X$ is smooth, is to try a line bundle that is ``handed'' to
us, namely some (positive or negative) power of the \emph{canonical line bundle};
$\omega_X=\det T_X^*$.  If $X$ is not smooth but instead normal, we can construct
$\omega_X$ on the smooth locus and then push it forward to obtain a rank one
reflexive sheaf on all of $X$ (which sometimes is still a line bundle).  Next we will
explore how we might ``force'' this line bundle to be ample in some (actually many)
cases.

Let $X$ be a minimal surface of general type that contains a
$(-2)$-curve (a smooth rational curve with self-intersection $-2$).
For an example of such a surface consider the following.
\begin{example}
  $\wtilde X=(x^5+y^5+z^5+w^5=0)\subseteq\bP^3$ with the $\bZ_2$-action that
  interchanges $x\leftrightarrow y$ and $z\leftrightarrow w$.  This action has five
  fixed points, $[1:1:-\varepsilon^i:-\varepsilon^i]$ for $i=1,\dots,5$ where
  $\varepsilon$ is a primitive $5^\text{th}$ root of unity.  Consequently the
  quotient $\factor\wtilde X.\bZ_2.$ has five singular points, each a simple double
  point of type $A_1$. Let $X\to \factor\wtilde X.\bZ_2.$ be the minimal resolution
  of singularities.  Then $X$ contains five $(-2)$-curves, the exceptional divisors
  over the singularities.
\end{example}

Let us return to the general case, that is, $X$ is a minimal surface
of general type that contains a $(-2)$-curve, $C\subseteq X$. As
$C\simeq\bP^1$, and $X$ is smooth, the adjunction formula gives us that
$K_X\cdot C=0$. Therefore $K_X$ is not ample.

On the other hand, since $X$ is a minimal surface of general type, it follows that
$K_X$ is semi-ample, that is, some multiple of it is base-point free. In other words,
there exists a morphism,
$$
|mK_X|: X\to X_\text{can}\subseteq\bP(\coh 0.X.\sO_X(mK_X).^{*}).
$$
This may be deduced from various results. For example, it follows from Bombieri's
classification of pluri-canonical maps, but perhaps the simplest proof is provided by
Miles Reid \cite[E.3]{MR1442522}.

It is then relatively easy to see that this morphism onto its image is independent of
$m$ (as long as $mK_X$ is base point free).  This constant image is called the
\emph{canonical model} of $X$, it will be denoted by $X_\text{can}$.

The good news is that the canonical line bundle of $X_\text{can}$ is indeed ample,
but the trouble is that $X_\text{can}$ is singular. We might consider this as the first
sign of the necessity of working with singular varieties. Fortunately the
singularities are not too bad, so we still have a good chance to work with this
model. In fact, the singularities that can occur on the canonical model of a surface
of general type belong to a much studied class. This class goes by several names;
they are called \emph{du Val singularities}, or \emph{rational double points}, or
\emph{Gorenstein, canonical singularities}. For more on these singularities, refer to
\cite{MR543555}, \cite{Reid87}.

\subsection{Normal crossings}
\label{subsectionNormalCrossings}

These singularities already appear in the construction of the moduli space of stable
curves (or if the reader prefers, the construction of a compactificaton of the moduli
space of smooth projective curves). If we want to understand degenerations of smooth
families, we have to allow normal crossings.

A \emph{normal crossing} singularity is one that is locally analytically (or
formally) isomorphic to the intersection of coordinate hyperplanes in a linear space.
In other words, it is a singularity locally analytically defined as $(x_1x_2\cdots
x_r=0)\subseteq\bA^n$ for some $r\leq n$. In particular, as opposed to the curve
case, for surfaces it allows for triple intersections.  However, triple intersections
may be ``resolved'':
Let $X=(xyz=0)\subseteq\bA^3$. Blow up the origin $O\in\bA^3$,
$\sigma:Bl_O\bA^3\to \bA^3$ and consider the proper transform of $X$,
$\sigma: \wtilde X\to X$. Observe that $\wtilde X$ has only double
normal crossings.

Another important point to remember about normal crossings is that
they are {\it not} normal. In particular they do not belong to the
previous category.
For some interesting and perhaps surprising examples of surfaces with
normal crossings see \cite{0705.0926v2}.

\subsection{Pinch points}

Another non-normal singularity that can occur as the limit of smooth varieties is the
pinch point. It is locally analytically defined as $(x_1^2=x_2 x_3^2)\subseteq\bA^n$.
This singularity is a double normal crossing away from the pinch point. Its
normalization is smooth, but blowing up the pinch point (i.e., the origin) does not
make it any better. (Try it for yourself!)

\subsection{Cones}

Let $C\subseteq \bP^2$ be a curve of degree $d$ and $X\subseteq \bP^3$
the projectivized cone over $C$. As $X$ is a degree $d$ hypersurface,
it admits a smoothing.

\begin{example}\label{ex:cone}
  Let $\Xi=(x^d+y^d+z^d+tw^d=0)\subseteq \bP^3_{x:y:z:w}\times
  \bA^1_t$. The special fiber $\Xi_0$ is a cone over a smooth plane
  curve of degree $d$ and the general fiber $\Xi_t$, for $t\neq 0$, is
  a smooth surface of degree $d$ in $\bP^3$.
\end{example}

This, again, suggests that we must allow some singularities. The question is, whether
we can limit the type of singularities we must deal with. More particularly to this
case, can we limit the type of cones we need to allow?

First we need an auxiliary computation. By the nature of the computation it is easier
to use \emph{divisors} instead of \emph{line bundles}.

\begin{demo}{Commentary}\label{comm:ampleness}
  One of our ultimate goals is to construct a moduli space for canonical models of
  varieties. We are already aware that the minimal model program has to deal with
  singularities and so we must allow some singularities on canonical models.  We
  would also like to understand what constraints are imposed if our goal is to
  construct a moduli space. The point is that in order to construct our moduli space,
  the objects must have an ample canonical class. It is possible that a family of
  canonical models degenerates to a singular fiber that has singularities worse than
  the original canonical models. An important question then is whether we may resolve
  the singularities of this special fiber and retain ampleness of the canonical
  class. The next example shows that this is not always possible.
\end{demo}

\begin{example}\label{ex:reducible-surface}
  Let $W$ be a smooth variety and $X=X_1\cup X_2\subseteq W$ such that $X_1$ and
  $X_2$ are Cartier divisors in $W$. Then by the adjunction formula we have
  $$
  \begin{aligned}
    K_X& = (K_W+X)\resto X \\
    K_{X_1} & = (K_W+X_1)\resto {X_1}\\
    K_{X_2} & = (K_W+X_2)\resto {X_2}
  \end{aligned}
  $$
  Therefore
  \begin{equation}
    \label{eq:3}
    K_X\resto{X_i}=K_{X_i}+X_{3-i}\resto{X_i}
  \end{equation}
  for $i=1,2$, so we have that
  \begin{equation}
    \label{eq:2}
    K_X \text{ is ample } \quad\Leftrightarrow\quad
    K_X\resto{X_i}=K_{X_i}+X_{3-i}\resto{X_i} \text{ is ample for }i=1,2.
  \end{equation}
\end{example}

Next, let $X$ be a normal projective surface with $K_X$ ample and an isolated
singular point $P\in\Sing X$. Assume that $X$ is isomorphic to a cone $\Xi_0\subseteq
\bP^3$ as in Example~\ref{ex:cone} locally analytically near $P$.  Further assume
that $X$ is the special fiber of a family $\Xi$ that itself is smooth. In particular,
we may assume that all fibers other than $X$ are smooth.  As explained in
\eqref{comm:ampleness}, we would like to see whether we may resolve the singular
point $P\in X$ and still be able to construct our desired moduli space, i.e., that
$K$ of the resolved fiber would remain ample. For this purpose we may assume that $P$
is the only singular point of $X$.

Let $\Upsilon\to \Xi$ be the blowing up of $P\in \Xi$ and let $\wtilde X$ denote the
proper transform of $X$. Then $\Upsilon_0=\wtilde X\cup E$ where $E\simeq \bP^2$ is
the exceptional divisor of the blow up.  Clearly, $\sigma:\wtilde X\to X$ is the blow
up of $P$ on $X$, so it is a smooth surface and $\wtilde X\cap E$ is isomorphic to
the degree $d$ curve over which $X$ is locally analytically a cone.

We would like to determine the condition on $d$ that ensures that the
canonical divisor of $\Upsilon_0$ is still ample. According to (\ref{eq:2})
this means that we need that $K_E+\wtilde X\resto E$ and $K_{\wtilde
  X}+E\resto{\wtilde X}$ be ample.

As $E\simeq \bP^2$, $\omega_E\simeq \sO_{\bP^2}(-3)$, so
$\sO_E(K_E+\wtilde X\resto E)\simeq \sO_{\bP^2}(d-3)$. This is ample if
and only if $d>3$.

As this computation is local near $P$ the only relevant issue about the ampleness of
$K_{\wtilde X}+E\resto{\wtilde X}$ is whether it is ample in a neighborhood of
$E_X\leteq E\resto{\wtilde X}$. By the next claim this is equivalent to asking when
$(K_{\wtilde X}+E_X)\cdot E_X$ is positive.

\begin{demo*}{\it Claim}
  Let $Z$ be a smooth projective surface with non-negative Kodaira dimension
  and $\Gamma\subset Z$ an effective divisor. If $(K_Z+\Gamma)\cdot
  C>0$ for every proper curve $C\subset Z$, then $K_Z+\Gamma$ is ample.
\end{demo*}

\begin{demo*}{\it Proof}
  By the assumption on the Kodaira dimension there exists an $m>0$ such that $mK_Z$
  is effective, hence so is $m(K_Z+\Gamma)$. Then by the assumption on the
  intersection number, $(K_Z+\Gamma)^2>0$, so the statement follows by the
  Nakai-Moishezon criterium.  \qed
\end{demo*}

Observe that by the adjunction formula $(K_{\wtilde X}+E_X)\cdot E_X=\deg K_{E_X} =
d(d-3)$ as $E_X$ is isomorphic to a plane curve of degree $d$.  Again, we obtain the
same condition as above and thus conclude that $K_{\Upsilon_0}$ may be ample only if
$d>3$.

Now, if we are interested in constructing moduli spaces, then one of the requirements
of being stable is that the canonical bundle be ample. This means that in order to
obtain a compact moduli space we have to allow cone singularities over curves of
degree $d\leq 3$. The singularity we obtain for $d=2$ is a rational double point, but
the singularity for $d=3$ is not even rational. This does not fit any of the earlier
classes we discussed. It belongs to the one discussed in the next section.

\subsection{Log canonical singularities}

Let us investigate the previous situation under more general
assumptions.

\begin{demo}{Computation}
  Let $D=\sum_{i=0}^r \lambda_iD_i$, ($\lambda_i\in\bN$), be a divisor with only
  normal crossing singularities in a smooth ambient variety such that $\lambda_0=1$.
  Using a generalized version of the adjunction formula shows that in this situation
  (\ref{eq:3}) remains true.
  \begin{equation}
    \label{eq:4}
    K_D\resto{D_0}=K_{D_0}+\sum_{i=1}^r \lambda_iD_{i}\resto{D_0}
  \end{equation}

  Let $f:\Xi\to B$ a projective family with $\dim B=1$, $\Xi$ smooth
  and $K_{\Xi_b}$ ample for all $b\in B$. Further let $X=\Xi_{b_0}$
  for some $b_0\in B$ a singular fiber and let $\sigma:\Upsilon\to\Xi$
  be an embedded resolution of $X\subseteq \Xi$. Finally let
  $Y=\sigma^*X=\wtilde X+\sum_{i=1}^r \lambda_iF_i$ where $\wtilde X$ is the
  proper transform of $X$ and $F_i$ are exceptional divisors for
  $\sigma$. We are interested in finding conditions that are necessary
  for $K_Y$ to remain ample.

  Let $E_i\leteq F_i\resto{\wtilde X}$ be the exceptional divisors for $\sigma:
  \wtilde X\to X$ and for the simplicity of computation, assume that the $E_i$ are
  irreducible.  For $K_Y$ to be ample we need that $K_Y\resto{\wt X}$ as well as
  $K_Y\resto{F_i}$ for all $i$ are all ample. Clearly, the important one of these for
  our purposes is $K_Y\resto{\wt X}$ for which by (\ref{eq:4}) we have that
  $$
  K_Y\resto{\wtilde X}=K_{\wtilde X}+\sum_{i=1}^r \lambda_iE_i.
  $$
  As usual, we may write $K_{\wtilde X}=\sigma^*K_X+\sum_{i=1}^r a_iE_i$, so we
  are looking for conditions to guarantee that $\sigma^*K_X+\sum
  (a_i+\lambda_i)E_i$ be ample. In particular, its restriction to any of the $E_i$
  has to be ample.
  To further simplify our computation let us assume that $\dim X=2$.  Then the
  condition that we want satisfied is that for all $j$,
  \begin{equation}
    \label{eq:5}
    \left(\sum_{i=1}^r(a_i+\lambda_i)E_i\right)\cdot E_j >0.
  \end{equation}

  Let
  $$
  \begin{aligned}
    E_+&=\sum_{a_i+\lambda_i\geq 0} \vert a_i+\lambda_i\vert E_i \text{,\quad and}\\
    E_-&=\sum_{a_i+\lambda_i< 0} \vert a_i+\lambda_i\vert E_i \text{,\quad so}\\
    \sum_{i=1}^r(a_i+\lambda_i)E_i&=E_+-E_-.
  \end{aligned}
  $$

  Choose a $j$ such that $E_j\subseteq \Supp E_+$. Then $E_-\cdot E_j\geq 0$ since
  $E_j\not\subseteq E_-$ and (\ref{eq:5}) implies that $(E_+-E_-)\cdot E_j>0$. These
  together imply that $E_+\cdot E_j>0$ and then that $E_+^2> 0$. However, the $E_i$
  are exceptional divisors of a birational morphism, so their intersection matrix,
  $(E_i\cdot E_j)$ is negative definite.

  The only way this can happen is if $E_+=0$. In other words, $a_i+\lambda_i< 0$ for
  all $i$.  However, the $\lambda_i$ are positive integers, so this implies that
  $K_Y$ may remain ample only if $a_i< -1$ for all $i=1,\dots,r$.

  The definition of a \emph{log canonical singularity} is the exact opposite of this
  condition. It requires that $X$ be normal and admit a resolution of singularities,
  say $Y\to X$, such that all the $a_i\geq -1$. This means that the above argument
  shows that we may stand a fighting chance if we resolve singularities that are
  \emph{worse} than log canonical, but have no hope to do so with log canonical
  singularities.  In other words, this is another class of singularities that we have
  to allow. As we remarked above, the class of singularities we obtained for the
  cones in the previous subsection belong to this class. In fact, all the normal
  singularities that we have considered so far belong to this class.

  The good news is that by now we have covered pretty much all the ways that
  something can go wrong and found the class of singularities we must allow. Since we
  have already found that we have to deal with some non-normal singularities and in
  fact in this example we have not really needed that $X$ be normal, we conclude that
  we will have to allow the non-normal cousins of log canonical singularities. These
  are called \emph{semi-log canonical singularities} and the reader can find their
  definition in the next subsection.
\end{demo}

\subsection{Semi-log canonical singularities}
\label{SubsectionSLC}

Semi-log canonical singularities are very important in moduli theory. These are
exactly the singularities that appear on stable varieties, the higher dimensional
analogs of stable curves. However, their definition is rather technical, so the
reader might want to skip this section at the first reading.

As a warm-up, let us first define the normal and more traditional
singularities that are relevant in the minimal model program.

\begin{defini}\label{def:semi-log-canonical}
  A pair $(X,\Delta)$ is called \emph{log~$\bQ$-Gorenstein} if $K_X+\Delta$ is
  $\bQ$-Cartier, i.e., some integer multiple of $K_X+\Delta$ is a Cartier divisor.
  Let $(X,\Delta)$ be a log~$\bQ$-Gorenstein pair and $f:\wtilde X\to X$ a log
  resolution of singularities with exceptional divisor $E=\union E_i$. Express the
  log canonical divisor of $\wtilde X$ in terms of $K_X+\Delta$ and the exceptional
  divisors:
  $$
  K_{\wtilde X} + \wt\Delta \equiv f^*(K_X +\Delta) + \sum a_i E_i
  $$
  where $\wt\Delta=f^{-1}_*\Delta$, the strict transform of $\Delta$ on $\wt X$ and
  $a_i\in \bQ$. Then the pair $(X,\Delta)$ has
  \begin{center}
    $\left.\begin{matrix}
      \text{\emph{terminal}}\\
      \text{\emph{canonical}}\\
      \text{\emph{plt}}\\
      \text{\emph{klt}}\\
      \text{\emph{log canonical}}
    \end{matrix} \right\}$
     singularities, if \quad
    $
    \begin{matrix}
      \text{for all log resolutions $f$,}\\
      \text{and for all $i$,}
    \end{matrix}
    $
    $\left\{\begin{matrix}
      a_i>0. \hfill \\
      a_i \geq 0. \hfill \\
      a_i> -1. \hfill \\
      a_i> -1 \text{ and } \rdown\Delta \leq 0.\\
      a_i\geq -1. \hfill
    \end{matrix}\right.$
  \end{center}
\end{defini}

The corresponding definitions for non-normal varieties are somewhat more cumbersome.
We include them here for completeness, but the reader should feel free to skip them
and assume that for instance ``semi-log canonical'' means something that can be
reasonably considered a non-normal version of log canonical.

Suppose that $X$ is a reduced equidimensional scheme that satisfies the following
conditions:
\begin{enumerate}
\item
  \label{item:1}
  $X$ satisfies Serre's condition S2 (cf.\ \cite[Thm.\ 8.22A(2)]{Hartshorne77}).
\item
  \label{item:2}
  $X$ has only simple normal double crossings in codimension 1 (in particular $X$ is
  Gorenstein in codimension 1)\footnote{Sometimes a ring that is S2 and Gorenstein in
    codimension 1 is called quasi-normal.}.
\end{enumerate}
Conditions (\ref{def:semi-log-canonical}.\ref{item:1}) and
(\ref{def:semi-log-canonical}.\ref{item:2}) imply that we may treat the canonical
module of $X$ as a divisorial sheaf even though $X$ is not normal.  Further suppose
that $D$ is a $\bQ$-Weil divisor on $X$ (again, following \cite[Chapter 16]{K+92}, we
assume that $X$ is regular at the generic point of each component in $\Supp D$).

\begin{remark}
  Notice that conditions (\ref{def:semi-log-canonical}.\ref{item:1}) and
  (\ref{def:semi-log-canonical}.\ref{item:2}) imply that $X$ is seminormal since it
  is seminormal in codimension 1 (see \cite[Corollary 2.7]{GrecoTraversoSeminormal}).
\end{remark}

Set $\rho : X^N \rightarrow X$ to be the normalization of $X$ and suppose that $B$ is
the divisor of the conductor ideal on $X^N$.  We use $\rho^{-1}(D)$ to denote the
pullback of $D$ to $X^N$.

\begin{definition}
  We say that $(X, D)$ is \emph{semi-log canonical} if the following two conditions
  hold.
  \begin{enumerate}
  \item
    $K_X + D$ is $\bQ$-Cartier, and
  \item
    the pair $(X^N, B + \rho^{-1} D)$ is log canonical.
  \end{enumerate}
\end{definition}

Actually, this is not the original definition of semi-log canonical singularities.
The original definition (which is equivalent to this one) uses the theory of semi-resolutions.  See
\cite{KollarShepherdBarron}, \cite[Chapter 12]{K+92}, and \cite{KollarSemiLogRes} for
details.

\section{Hyperresolutions and Du~Bois' original definition}

A very important construction is Du~Bois's generalized De Rham complex. The original
construction of Du~Bois's complex, $\FullDuBois{X}$, is based on simplicial
resolutions. The reader interested in the details is referred to the original article
\cite{DuBoisMain}.  Note also that a simplified construction was later obtained in
\cite{Carlson85} and \cite{GNPP88} via the general theory of polyhedral and cubic
resolutions.  At the end of the paper, we include an appendix in which we explain how
to construct, and give examples of cubical hyperresolutions.  An easily accessible
introduction can be found in \cite{SteenbrinkgVanishing}.  Another useful reference
is the recent book \cite{PetersSteenbrinkBook}.

Recently the second named author found a simpler alternative construction of (part
of) Du~Bois's complex that does not need a simplicial resolution, see
\cite{SchwedeEasyCharacterization} and also Section
\ref{SectionHyperresolutionFreeCharacterizations} below.  However we will discuss the
original construction because we believe that it is important to keep in mind the way
these singularities appeared as that explains their usefulness.  For more on
applications of Du~Bois's complex and Du~Bois singularities see
\cite{SteenbrinkMixed}, \cite[Chapter 12]{KollarShafarevich},
\cite{KovacsDuBoisLC1,KovacsDuBoisLC2}.

The word ``hyperresolution'' will refer to either simplicial, polyhedral, or cubic
resolution. Formally, the construction of $\FullDuBois{X}$ is the same regardless the
type of resolution used and no specific aspects of either types will be used.

The following definition is included to make sense of the statements of some of the
forthcoming theorems. It can be safely ignored if the reader is not interested in the
detailed properties of Du~Bois's complex and is willing to accept that it is a very
close analog of the De~Rham complex of smooth varieties.

\begin{defini}
  Let $X$ be a complex scheme (i.e., a scheme of finite type over $\bC$) of dimension
  n. Let $D_{\rm filt}(X)$ denote the derived category of filtered complexes of
  $\sO_{X}$-modules with differentials of order $\leq 1$ and $D_{filt, coh}(X)$ the
  subcategory of $D_{\rm filt}(X)$ of complexes $\cx K$, such that for all $i$, the
  cohomology sheaves of $Gr^{i}_{\rm filt}K^{\mydot}$ are coherent cf.\
  \cite{DuBoisMain}, \cite{GNPP88}.  Let $D(X)$ and $D_{\rm coh}(X)$ denote the
  derived categories with the same definition except that the complexes are assumed
  to have the trivial filtration.  The superscripts $+, -, b$ carry the usual meaning
  (bounded below, bounded above, bounded).  Isomorphism in these categories is
  denoted by $\qis$.  A sheaf $\sF$ is also considered a complex $\sF^\mydot$ with
  $\sF^0=\sF$ and $\sF^i=0$ for $i\neq 0$.  If $K^{\mydot}$ is a complex in any of
  the above categories, then $h^i(K^{\mydot})$ denotes the $i$-th cohomology sheaf of
  $K^{\mydot}$.

  The right derived functor of an additive functor $F$, if it exists, is denoted by
  $RF$ and $R^iF$ is short for $h^i\circ RF$. Furthermore, $\bH^i$, $\bH^i_Z$ , and
  $\sH^i_Z$ will denote $R^i\Gamma$, $R^i\Gamma_Z$, and $R^i\sH_Z$ respectively,
  where $\Gamma$ is the functor of global sections, $\Gamma_Z$ is the functor of
  global sections with support in the closed subset $Z$, and $\sH_Z$ is the functor
  of the sheaf of local sections with support in the closed subset $Z$. Note that
  according to this terminology, if $\phi\col Y\to X$ is a morphism and $\sF$ is a
  coherent sheaf on $Y$, then $R\phi_*\sF$ is the complex whose cohomology sheaves
  give rise to the usual higher direct images of $\sF$.
\end{defini}

\begin{thm}[{\cite[6.3, 6.5]{DuBoisMain}}]\label{defDB}
  Let $X$ be a proper complex scheme of finite type and $D$ a closed subscheme whose
  complement is dense in $X$. Then there exists a unique object $\FullDuBois X \in
  \Ob D_{\rm filt}(X)$ such that using the notation
  $$
  \Om^p_X\leteq Gr^{p}_{\rm filt}\, \FullDuBois{X}[p],
  $$
  it satisfies the following properties
  \begin{enumerate}
  \item 
    $ \FullDuBois{X} \qis \bC_{X} $, i.e., $\FullDuBois{X}$ is a resolution of the
    constant sheaf $\bC$ on $X$.

  \item $\underline{\Omega}_{(\blank)}^{\mydot}$ is functorial, i.e., if $\phi \col
    Y\to X$ is a morphism of proper complex schemes of finite type, then there exists
    a natural map $\phi^{*}$ of filtered complexes
    $$
    \phi^{*}\col \FullDuBois{X} \to R\phi_{*}\underline{\Omega}_Y^{\mydot}.
    $$
    Furthermore, $\FullDuBois{X} \in \Ob \left(D^{b}_{filt, coh}(X)\right)$ and if
    $\phi$ is proper, then $\phi^{*}$ is a morphism in $D^{b}_{filt, coh}(X)$.
    \label{functorial}

  \item Let $U \subseteq X$ be an open subscheme of $X$. Then
    $$
    \FullDuBois{X}\resto U 
    \qis\underline{\Omega}^{\,\mydot}_U.
    $$

  \item If $X$ is proper, there exists a spectral sequence degenerating at $E_1$ and abutting to the
    singular cohomology of $X$:
    $$
    E_1^{pq}={\bH}^q \left(X, \Om^p_X\right) \Rightarrow H^{p+q}(X^{\an}, \bC).
    $$\label{item:Hodge}

  \item If\/ $\varepsilon_\mydot\col X_\mydot\to X$ is a hyperresolution, then
    $$
    \FullDuBois{X}\qis \myR{\varepsilon_\mydot}_* \Omega^\mydot_{X_\mydot}.
    $$
    In particular, $h^i\left(\Om^p_X\right)=0$ for $i<0$.

  \item There exists a natural map, $\sO_{X}\to \Om^0_X$, compatible with
    $(\ref{defDB}.\ref{functorial})$. \label{item:dR-to-DB}

  \item If\/ $X$ is smooth, then
    $$
    \FullDuBois{X}\qis\Omega^\mydot_X.
    $$
    In particular,
    $$
    \Om^p_X\qis\Omega^p_X.
    $$

  \item If\/ $\phi\col Y\to X$ is a resolution of singularities, then
    $$
    \Om_X^{\dim X}\qis \myR\phi_*\omega_Y.
    $$
  \item\label{item:exact-triangle}
    Suppose that $\pi : \tld Y \rightarrow Y$ is a projective morphism and $X \subset
    Y$ a reduced closed subscheme such that $\pi$ is an isomorphism outside of $X$.
    Let $E$ denote the reduced subscheme of $\tld Y$ with support equal to
    $\pi^{-1}(X)$ and $\pi' : E \rightarrow X$ the induced map.  Then for each $p$
    one has an exact triangle of objects in the derived category,
    $$
    \xymatrix{ \Om^p_Y \ar[r] & \Om^p_X \oplus \myR \pi_* \Om^p_{\tld Y} \ar[r]^-{-}
      & \myR
      \pi'_* \Om^p_E \ar[r]^-{+1} & .\\
    }
    $$
  \end{enumerate}
\end{thm}

It turns out that Du~Bois's complex behaves very much like the de~Rham complex for
smooth varieties. Observe that (\ref{defDB}.\ref{item:Hodge}) says that the
Hodge-to-de~Rham spectral sequence works for singular varieties if one uses the
Du~Bois complex in place of the de~Rham complex. This has far reaching consequences
and if the associated graded pieces, $\Om^p_X$ turn out to be computable, then this
single property leads to many applications.

Notice that (\ref{defDB}.\ref{item:dR-to-DB}) gives a natural map $\sO_{X}\to
\Om^0_X$, and we will be interested in situations when this map is a
quasi-isomorphism.  When $X$ is proper over $\bC$, such a quasi-isomorphism will
imply that the natural map
\[
H^i(X^{\an}, \bC) \rightarrow H^i(X, \O_{X}) = \bH^i(X, \DuBois{X})
\]
is surjective because of the degeneration at $E_1$ of the spectral sequence in
(\ref{defDB}.\ref{item:Hodge}).

Following Du~Bois, Steenbrink was the first to study this condition and he christened
this property after Du~Bois.

\begin{defini}
  A scheme $X$ is said to have \emph{Du~Bois singularities} (or \emph{DB
    singularities} for short) if the natural map $\sO_{X}\to \Om^0_X$ from
  (\ref{defDB}.\ref{item:dR-to-DB}) is a quasi-isomorphism.
\end{defini}

\begin{remark}
  If $\varepsilon : X_{\mydot} \rightarrow X$ is a hyperresolution of $X$ (see the
  Appendix for a how to construct cubical hyperresolutions) then $X$ has Du~Bois
  singularities if and only if the natural map $\sO_X \rightarrow \myR
  {\varepsilon_{\mydot}}_* \sO_{X_{\mydot}}$ is a quasi-isomorphism.
\end{remark}

\begin{example}
  It is easy to see that smooth points are Du~Bois and Deligne proved that normal
  crossing singularities are Du~Bois as well cf.\ \cite[Lemme 2(b)]{MR0376678}.
\end{example}

We will see more examples of Du~Bois singularities in later sections.

\section{An injectivity theorem and splitting the Du~Bois complex}

In this section, we state an injectivity theorem involving the dualizing sheaf
that plays a role for Du~Bois singularities similar to the role that
Grauert-Riemenschneider players for rational singularities.  As an application, we
state a criterion for Du~Bois singularities related to a ``splitting'' of the Du~Bois
complex.

\begin{theorem}\protect{\cite[Lemma 2.2]{KovacsDuBoisLC1}, \cite[Proposition
  5.11]{SchwedeFInjectiveAreDuBois}}\label{TheoremKeyInjectivity} %
  Let $X$ be a reduced scheme of finite type over $\bC$, $x\in X$ a (possibly
  non-closed) point, and $Z=\overline{\{x\}}$ its closure. Assume that %
  $X\setminus Z$
  %
  has Du~Bois singularities in a neighborhood of $x$ (for example, $x$ may correspond
  to an irreducible component of the non-Du~Bois locus of $X$).
  Then the natural map
  \[
  \myH^i \left(\myR \sHom_X^{\mydot}(\DuBois{X}, \omega_X^{\mydot}) \right)_x
  \rightarrow \myH^i \left( \omega_X^{\mydot} \right)_x
  \]
  is injective for every $i$.
\end{theorem}

The proof uses the fact that for a projective $X$, $H^i(X^{\an}, \bC) \rightarrow
\bH^i(X, \DuBois{X})$ is surjective for every $i > 0$, which follows from Theorem
\ref{defDB}.

It would also be interesting and useful if the following generalization of this
injectivity were true.
\begin{question}\label{qu:injective}
  Suppose that $X$ is a reduced scheme essentially of finite type over $\bC$.  Is it
  true that the natural map of sheaves
  \[
  \myH^i \left(\myR \sHom_X^{\mydot}(\DuBois{X}, \omega_X^{\mydot}) \right)
  \rightarrow \myH^i \left( \omega_X^{\mydot} \right)
  \]
  is injective for every $i$?
\end{question}

Even though Theorem \ref{TheoremKeyInjectivity} does not answer
Question~\ref{qu:injective}, it has the following extremely useful corollary.

\begin{theorem} \cite[Theorem 2.3]{KovacsDuBoisLC1}, \cite[Theorem
  12.8]{KollarShafarevich} \label{ThmDuBoisSplitting}%
  Suppose that the natural map $\sO_X \rightarrow \DuBois{X}$ has a left inverse in
  the derived category (that is, a map $\rho : \DuBois{X} \rightarrow \sO_X$ such
  that the composition $\xymatrix{\sO_X \ar[r] & \DuBois{X} \ar[r]^{\rho} & \sO_X}$
  is an isomorphism).  Then $X$ has Du~Bois singularities.
\end{theorem}
\begin{proof}
  Apply the functor $\myR \sHom_{X}(\blank, \omega_X^{\mydot})$ to the maps
  $\xymatrix{\sO_X \ar[r] & \DuBois{X} \ar[r]^{\rho} & \sO_X}$.  Then by the
  assumption,
  the composition
  \[
  \xymatrix{ \omega_{X}^{\mydot} \ar[r]^-{\delta} & \myR \sHom_X(\DuBois{X},
    \omega_{X}^{\mydot}) \ar[r] & \omega_{X}^{\mydot} }
  \]
  is an isomorphism. Let $x\in X$ be a possibly non-closed point corresponding to an
  irreducible component of the non-Du~Bois locus of $X$ and consider the stalks at
  $x$ of the cohomology sheaves of the complexes above. We obtain that the natural
  map
  \[
  \myH^i \left(\myR \sHom_X(\DuBois{X}, \omega_X^{\mydot}) \right)_x
  \rightarrow \myH^i \left( \omega_X^{\mydot} \right)_x
  \]
  is surjective for every $i$.  But it is also injective by Theorem
  \ref{TheoremKeyInjectivity}. This proves that
  $\delta: (\omega_{X}^{\mydot})_x \to \myR \sHom_X(\DuBois{X}, \omega_{X}^{\mydot})_x$ is
  a quasi-isomorphism.  Finally, applying the functor $\myR \sHom_{\O_{X, x}}(\blank,
  (\omega_X^{\mydot})_x)$ one more time proves that $X$ is Du~Bois at $x$, contradicting our choice of $x \in X$
\end{proof}

This also gives the following Boutot-like theorem for Du~Bois singularities (cf.\
\cite{BoutotRational}).

\begin{corollary} \cite[Theorem 2.3]{KovacsDuBoisLC1}, \cite[Theorem
  12.8]{KollarShafarevich}\label{BoutoutDuBois} Suppose that $f : Y \rightarrow X$ is
  a morphism, $Y$ has Du~Bois singularities and the natural map $\sO_X \rightarrow
  \myR f_* \sO_Y$ has a left inverse in the derived category.  Then $X$ also has
  Du~Bois singularities.
\end{corollary}
\begin{proof}
  Observe that the composition is an isomorphism
  \[
  \sO_X \rightarrow \DuBois{X} \rightarrow Rf_*\DuBois{Y} \simeq \myR f_* \sO_Y
  \rightarrow \sO_X.
  \]
  Then apply Theorem \ref{ThmDuBoisSplitting}.
\end{proof}

As an easy corollary, we see that rational singularities are Du~Bois (which was first
observed in the isolated case by Steenbrink in \cite[Proposition
3.7]{SteenbrinkMixed}).

\begin{corollary}\cite{KovacsDuBoisLC1}, \cite{SaitoMixedHodge}
  If $X$ has rational singularities, then $X$ has Du~Bois singularities.
\end{corollary}
\begin{proof}
  Let $\pi : \tld X \rightarrow X$ be a log resolution.  One has the following
  composition $\sO_X \rightarrow \DuBois{X} \rightarrow \myR \pi_* \sO_{\tld X}$.
  Since $X$ has rational singularities, this composition is a quasi-isomorphism.
  Apply Corollary \ref{BoutoutDuBois}.
\end{proof}


\section{
Hyperresolution-free characterizations of Du~Bois singularities}
\label{SectionHyperresolutionFreeCharacterizations}

The definition of Du~Bois singularities given via hyperresolution is relatively
complicated (hyperresolutions themselves can be rather complicated to compute, see
\ref{AppendixOnHyper}).  In this section we state several hyperresolution free
characterizations of Du~Bois singularities.  The first such characterization was
given by Steenbrink in the isolated case. Another, more analytic characterization was
given by Ishii and improved by Watanabe in the isolated quasi-Gorenstein\footnote{A
  variety $X$ is quasi-Gorenstein if $K_X$ is a Cartier divisor.  It is not required
  that $X$ is Cohen-Macaulay.} case. Finally the second named author gave a
characterization that works for any reduced scheme.

Du~Bois gave a relatively simple characterization of an affine cone over a projective
variety being Du~Bois in \cite{DuBoisMain}.  Steenbrink generalized this criterion to
all normal isolated singularities.  It is this criterion that Steenbrink, Ishii,
Watanabe, and others used extensively to study isolated Du~Bois singularities.

\begin{theorem}\cite[Proposition 4.16]{DuBoisMain}
  \cite[3.6]{SteenbrinkMixed}
  \label{TheoremSteenbrinkIsolatedCharacterization}
  Let $(X, x)$ be a normal isolated Du~Bois singularity, and $\pi : \tld X
  \rightarrow X$ a log resolution of $(X, x)$ such that $\pi$ is an isomorphism
  outside of $X \setminus \{x\}$. Let $E$ denote the reduced pre-image of $x$.  Then
  $(X, x)$ is a Du~Bois singularity if and only if the natural map
  \[
  \myR^i \pi_* \sO_{\tld X} \rightarrow \myR^i \pi_* \sO_E
  \]
  is an isomorphism for all $i > 0$.
\end{theorem}
\begin{proof}
  Using Theorem \ref{defDB}, we have an exact triangle
  \[
  \xymatrix{ \DuBois{X} \ar[r] & \DuBois{\{x\}} \oplus \myR \pi_* \DuBois{\tld X}
    \ar[r]^-{-}
    & \myR \pi_* \DuBois{E} \ar[r]^-{+1} & .\\
  }
  \]
  Since $\{x\}$, $\tld X$ and $E$ are all Du~Bois (the first two are smooth, and $E$
  is snc), we have the following exact triangle
  \[
  \xymatrix{ \DuBois{X} \ar[r] & \sO_{\{x\}} \oplus \myR \pi_* \sO_{\tld X}
    \ar[r]^-{-} & \myR
    \pi_* \sO_{E} \ar[r]^-{+1} & .\\
  }
  \]
  Suppose first that $X$ has Du~Bois singularities (that is, $\DuBois{X} \qis \sO_X$).
  By taking cohomology and examining the long exact sequence, we see that $\myR^i
  \pi_* \sO_{\tld X} \rightarrow \myR^i \pi_* \sO_E$ is an isomorphism for all $i >
  0$.

  So now suppose that $\myR^i \pi_* \sO_{\tld X} \rightarrow \myR^i \pi_* \sO_E$ is
  an isomorphism for all $i > 0$.  By considering the long exact sequence of
  cohomology, we see that $\myH^i(\DuBois{X}) = 0$ for all $i > 0$.  On the other
  hand $\myH^0(\DuBois{X})$ is naturally identified with the seminormalization of
  $\sO_X$, see Proposition \ref{SeminormalDuBoisComplex} below.  Thus if $X$ is
  normal, then $\sO_X \rightarrow \myH^0(\DuBois{X})$ is an isomorphism.
\end{proof}

We now state a more analytic characterization, due to Ishii and slightly improved by
Watanabe.  First we recall the definition of the plurigenera of a singularity.
\begin{definition}
  For a singularity $(X, x)$, we define the plurigenera $\{\delta_m\}_{m \in \bN}$;
  \[
  \delta_m(X, x) = \dim_{\bC}\Gamma(X \setminus {x}, \O_X(mK_X)) / L^{2/m}( X
  \setminus \{x\}),
  \]
  where $L^{2/m}( X \setminus \{x\})$ denotes the set of all $L^{2/m}$-integrable
  $m$-uple holomorphic $n$-forms on $X \setminus \{x\}$.
\end{definition}

\begin{theorem}\cite[Theorem 2.3]{Ishii85} \cite[Theorem 4.2]{WatanabeKimioPlurigenera}
  Let $f: \tld X \rightarrow X$ be a log resolution of a normal isolated Gorenstein
  singularity $(X, x)$ of dimension $n \geq 2$. Set $E$ to be the reduced exceptional
  divisor (the pre-image of $x$).  Then $(X, x)$ is a Du~Bois singularity if and only
  if $\delta_m(X,x) \leq 1$ for any $m \in \bN$.
\end{theorem}

In \cite{SchwedeEasyCharacterization}, the second named author gave a
characterization of arbitrary Du~Bois singularities that did not rely on
hyperresolutions, but instead used a single resolution of singularities.  An
improvement of this was also obtained in \cite[Proposition
2.20]{SchwedeTakagiRationalPairs}.  We provide a proof for the convenience of the
reader.

\begin{theorem} \cite{SchwedeEasyCharacterization}, \cite[Proposition
  2.20]{SchwedeTakagiRationalPairs}
  \label{EasyDuBoisCriterion}
  Let $X$ be a reduced separated scheme of finite type over a field of characteristic
  zero.  Suppose that $X \subseteq Y$ where $Y$ is smooth and suppose that $\pi :
  \tld Y \rightarrow Y$ is a proper birational map with $\tld Y$ smooth and where
  $\overline X =\pi^{-1}(X)_{\red}$, the reduced pre-image of $X$, is a simple normal
  crossings divisor (or in fact any scheme with Du~Bois singularities).  Then $X$ has
  Du~Bois singularities if and only if the natural map $\sO_X \rightarrow \myR \pi_*
  \sO_{\overline X}$ is a quasi-isomorphism.

  In fact, we can say more.  There is an isomorphism $\xymatrix{\myR \pi_*
    \sO_{\overline X} \ar[r]^{\sim} & \DuBois{X} }$ such that the natural map $\sO_X
  \rightarrow \DuBois{X}$ can be identified with the natural map $\sO_X \rightarrow
  \myR \pi_* \sO_{\overline X}$.
\end{theorem}
\begin{proof}
  We first assume that $\pi$ is an isomorphism outside of $X$.  Then using Theorem
  \ref{defDB}, we have an exact triangle
  \[
  \xymatrix{ \DuBois{Y} \ar[r] & \DuBois{X} \oplus \myR \pi_* \DuBois{\tld Y}
    \ar[r]^-{-} &
    \myR \pi_* \DuBois{\overline{X}} \ar[r]^-{+1} & .\\
  }
  \]
  Using the octahedral axiom, we obtain the following diagram
  \[
  \xymatrix{ C^{\mydot} \ar[d]^{\sim} \ar[r] & \DuBois{Y} \ar[d]^{\alpha} \ar[r] &
    \DuBois{X} \ar[d]^{\beta} \ar[r]^-{+1} & \\ C^{\mydot} \ar[r] & \myR \pi_*
    \DuBois{\tld Y} \ar[r] & \myR \pi_* \DuBois{\overline{X}} \ar[r]^-{+1} & .  }
  \]
  where $C^{\mydot}$ is simply the object in the derived category that completes the
  triangles.  But notice that the vertical arrow $\alpha$ is an isomorphism since $Y$
  has rational singularities (in which case each term in the middle column is
  isomorphic to $\sO_Y$).  Thus the vertical arrow $\beta$ is also an isomorphism.

  One always has a commutative diagram (where the arrows are the natural ones)
  \[
  \xymatrix{
    \sO_X \ar[d] \ar[r] & \DuBois{X} \ar[d]^{\beta} \\
    \myR \pi_* \sO_{\overline{X}} \ar[r]_{\delta} & \myR \pi_* \DuBois{\overline X} }
  \]
  Observe that $\overline{X}$ has Du~Bois singularities since it has normal
  crossings, thus $\delta$ is a quasi-isomorphism.  But then the theorem is proven at
  least in the case that $\pi$ is an isomorphism outside of $X$.

  For the general case, it is sufficient to show that $\myR \pi_* \sO_{\overline X}$
  is independent of the choice of resolution.  Since any two log resolutions can be
  dominated by a third, it is sufficient to consider two log resolutions $\pi_1 : Y_1
  \rightarrow Y$ and $\pi_2 : Y_2 \rightarrow Y$ and a map between them $\rho : Y_2
  \rightarrow Y_1$ over $Y$.  Let $F_1 = (\pi_1^{-1}(X))_{\red}$ and $F_2 =
  (\pi_2^{-1}(X))_{\red} = (\rho^{-1}(F_1))_{\red}$.  Dualizing the map and applying
  Grothendieck duality implies that it is sufficient to prove that $\omega_{Y_1}(F_1)
  \leftarrow \myR \rho_* ( \omega_{Y_2}(F_2))$ is a quasi-isomorphism.

  We now apply the projection formula while twisting by $\omega_{Y_1}^{-1}(-F_1)$.
  Thus it is sufficient to prove that
  \[
  \myR \rho_* (\omega_{Y_2/Y_1}(F_2 - \rho^* F_1)) \rightarrow \sO_{Y_1}
  \]
  is a quasi-isomorphism.  But note that $F_2 - \rho^* F_1 = - \lfloor \rho^*
  (1-\varepsilon) F_1 \rfloor$ for sufficiently small $\varepsilon > 0$.  Thus it is
  sufficient to prove that the pair $(Y_1, (1-\varepsilon)F_1)$ has klt singularities
  by Kawamata-Viehweg vanishing in the form of local vanishing for multiplier ideals;
  see \cite[9.4]{LazarsfeldPositivity2}.  But this is true since $Y_1$ is smooth and
  $F_1$ is a reduced integral divisor with simple normal
  crossings.
\end{proof}

It seems that in this characterization the condition that the ambient variety $Y$ is
smooth is asking for too much. We propose that the following may be a more natural
characterization. For some motivation and for a statement that may be viewed as a
sort of converse, see Conjecture~\ref{conj:DB-to-rtl} and the discussion preceding
it.
\begin{conjecture}
  Theorem \ref{EasyDuBoisCriterion} should remain true if the hypothesis that $Y$ is
  smooth is replaced by the condition that $Y$ has rational singularities.
\end{conjecture}

Having Du~Bois singularities is a local condition, so even if $X$ is not embeddable
in a smooth scheme, one can still use Theorem~\ref{EasyDuBoisCriterion} by passing to
an affine open covering.

To illustrate the utility and meaning of Theorem \ref{EasyDuBoisCriterion}, we will
explore the situation when $X$ is a hypersurface inside a smooth scheme $Y$.  Using
the notation of Theorem \ref{EasyDuBoisCriterion}, we note that we have the following
diagram of exact triangles.
\[
\xymatrix{ & \myR \pi_* \sO_{\tld Y}(-\overline{X}) \ar[r] & \myR \pi_* \sO_{\tld Y}
  \ar[r]
  & \myR \pi_* \sO_{\overline X} \ar[r]^-{+1} & \\
  0 \ar[r] & \sO_Y(-X) \ar[u]_{\alpha} \ar[r] & \sO_Y \ar[u]_{\beta} \ar[r] & \sO_X
  \ar[u]_{\gamma} \ar[r] & 0 }
\]
Since $Y$ is smooth, $\beta$ is a quasi-isomorphism (as then $Y$ has at worst
rational singularities).  Therefore, $X$ has Du~Bois singularities if and only if the
map $\alpha$ is a quasi-isomorphism.  However, $\alpha$ is a quasi-isomorphism if and
only if the dual map
\begin{equation}
  \label{EquationDualMap}
  \myR \pi_* \omega_{\tld Y}^{\mydot}(\overline X) \rightarrow \omega_Y^{\mydot}(X)
\end{equation}
is a quasi-isomorphism.  The projection formula tells us that Equation
\ref{EquationDualMap} is a quasi-isomorphism if and only
\begin{equation}
  \label{EquationDualMap2}
  \myR \pi_* \sO_{\tld Y}(K_{\tld Y/ Y} - \pi^* X + \overline X) \rightarrow \sO_X
\end{equation}
is a quasi-isomorphism.  Note however that $-\pi^* X + \overline X = \lceil -(1 -
\varepsilon) \pi^* X \rceil$ for $\varepsilon > 0$ and sufficiently close to zero.
Thus the left side of Equation \ref{EquationDualMap2} can be viewed as $\myR \pi_*
\sO_{\tld Y}(\lceil K_{\tld Y/ Y} -(1 - \varepsilon) \pi^* X \rceil)$ for
$\varepsilon > 0$ sufficiently small.  Note that Kawamata-Viehweg vanishing in the
form of local vanishing for multiplier ideals implies that $\mJ(Y, (1-\varepsilon)X)
\qis \myR \pi_* \sO_{\tld Y}(\lceil K_{\tld Y/ Y} -(1 - \varepsilon) \pi^* X
\rceil)$.  Therefore $X$ has Du~Bois singularities if and only if $\mJ(Y,
(1-\varepsilon)X) \simeq \sO_X$.

\begin{corollary}
  If $X$ is a hypersurface in a smooth $Y$, then $X$ has Du~Bois singularities if and
  only if the pair $(Y, X)$ is log canonical.
\end{corollary}

Note that Du~Bois hypersurfaces have also been characterized via the Bernstein-Sato
polynomial, see \cite[Theorem 0.5]{SaitoOnTheHodgeFiltrationOfHodgeModules}.

\section{Seminormality of Du~Bois singularities}

In this section we show that Du~Bois singularities are partially characterized by
seminormality.  First we remind the reader what it means for a scheme to be
seminormal.

\begin{definition} \cite{SwanSeminormality} \cite{GrecoTraversoSeminormal} Suppose
  that $R$ is a reduced excellent ring and that $S \supseteq R$ is a reduced
  $R$-algebra which is finite as an $R$-module.  We say that the extension $i : R
  \into S$ is \emph{subintegral} if the following two conditions are satisfied.
\begin{enumerate}
\item
  $i$ induces a bijection on spectra, $\Spec S \rightarrow \Spec R$.
\item
  $i$ induces an isomorphism of residue fields over every (not necessarily
  closed) point of $\Spec R$.
\end{enumerate}
\end{definition}

\begin{remark}
  In \cite{GrecoTraversoSeminormal}, subintegral extensions are called
  quasi-isomorphisms.
\end{remark}

\begin{definition} \cite{SwanSeminormality} \cite{GrecoTraversoSeminormal} Suppose
  that $R$ is a reduced excellent ring.  We say that $R$ is \emph{seminormal} if
  every subintegral extension $R \into S$ is an isomorphism.  We say that a scheme
  $X$ is \emph{seminormal} if all of its local rings are seminormal.
\end{definition}

\begin{remark}
  In \cite{GrecoTraversoSeminormal}, the authors call $R$ seminormal if there is no
  proper subintegral extension $R\into S$ such that $S$ is contained in the integral
  closure of $R$ (in its total field of fractions).  However, it follows from
  \cite[Corollary 3.4]{SwanSeminormality} that the above definition is equivalent.
\end{remark}

\begin{remark}
  Seminormality is a local property.  In particular, a ring is seminormal if and only
  if it is seminormal after localization at each of its prime (equivalently, maximal)
  ideals.
\end{remark}

\begin{remark}
  The easiest example of seminormal schemes are schemes with snc singularities.  In
  fact, a one dimensional variety over an algebraically closed field is seminormal if
  and only if its singularities are locally analytically isomorphic to a union of coordinate
  axes in affine space.
\end{remark}

We will use the following well known fact about seminormality.

\begin{lemma}
\label{LemmaSectionsAreSeminormal}
If $X$ is a seminormal scheme and $U \subseteq X$ is any open set, then $\Gamma(U,
\sO_X)$ is a seminormal ring.
\end{lemma}
\begin{proof}  We leave it as an exercise to the reader. \end{proof}

It is relatively easy to see, using the original definition via hyperresolutions,
that if $X$ has Du~Bois singularities, then it is seminormal.  Du~Bois certainly knew
this fact, see \cite[Proposition 4.9]{DuBoisMain} although he didn't use the word
seminormal.  Later Saito \cite{SaitoMixedHodge} proved that seminormality in fact
partially characterizes Du~Bois singularities.  We give a different proof of this
fact, due to the second named author.

\begin{proposition}\cite[Proposition 5.2]{SaitoMixedHodge} \cite[Lemma
  5.6]{SchwedeFInjectiveAreDuBois}
  \label{SeminormalDuBoisComplex}
  Suppose that $X$ is a reduced separated scheme of finite type over $\bC$.  Then
  $h^0(\DuBois{X}) = \sO_{X^{\sn}}$ where $\sO_{X^{\sn}}$ is the structure sheaf of
  the seminormalization of $X$.
\end{proposition}
\begin{proof}
  Without loss of generality we may assume that $X$ is affine.  We need only consider
  $\pi_* \sO_E$ by Theorem \ref{EasyDuBoisCriterion}.  By Lemma
  \ref{LemmaSectionsAreSeminormal}, $\pi_* \sO_E$ is a sheaf of seminormal rings.
  Now let $X' = \Spec (\pi_* \sO_E)$ and consider the factorization
  \[
  E \rightarrow X' \rightarrow X.
  \]
  Note $E \rightarrow X'$ must be surjective since it is dominant by construction and
  is proper by \cite[II.4.8(e)]{Hartshorne77}.  Since the composition has connected
  fibers, so must have $\rho : X' \rightarrow X$.  On the other hand, $\rho$ is a
  finite map since $\pi$ is proper.  Therefore $\rho$ is a bijection on points.
  Because these maps and schemes are of finite type over an algebraically closed
  field of characteristic zero, we see that $\Gamma(X, \sO_X) \rightarrow \Gamma(X',
  \sO_{X'})$ is a subintegral extension of rings.  Since $X'$ is seminormal, so is
  $\Gamma(X', \sO_{X'})$, which completes the proof.
\end{proof}

\section{A multiplier-ideal-like characterization of Cohen-Macaulay Du~Bois
  singularities}
\label{SectionAMultiplierIdealLikeCharacterization}
In this section we state a characterization of Cohen-Macaulay Du~Bois singularities
that explains why Du~Bois singularities are so closely linked to rational and log
canonical singularities.

We first do a suggestive computation.  Suppose that $X$ embeds into a smooth scheme
$Y$ and that $\pi : \tld Y \rightarrow Y$ is an embedded resolution of $X$ in $Y$
that is an isomorphism outside of $X$.  Set $\tld X$ to be the strict transform of
$X$ and set $\overline X$ to be the reduced pre-image of $X$.  We further assume that
$\overline X = \tld X \cup E$ where $E$ is a reduced simple normal crossings divisor
that intersects $\tld X$ transversally in another reduced simple normal crossing
divisor.  Note that $E$ is the exceptional divisor of $\pi$ (with reduced scheme
structure).  Set $\Sigma \subseteq X$ be the image of $E$.  We have the following
short exact sequence.
\[
0 \rightarrow \sO_{\tld X}(-E) \rightarrow \sO_{\overline X} \rightarrow \sO_{E}
\rightarrow 0
\]
We apply $\myR \sHom_{\sO_Y}(\blank, \omega_{\tld Y}^{\mydot})$ followed by $\myR
\pi_* \blank $ and obtain the following exact triangle.
\[
\xymatrix{ \myR \pi_* \omega_E^{\mydot} \ar[r] & \myR \pi_* \omega_{\overline
    X}^{\mydot} \ar[r] & \myR \pi_* \omega_{\tld X}(E)[\dim X] \ar[r]^-{+1}
  & }
\]
Using (\ref{defDB}.\ref{item:exact-triangle}), the left-most object can be identified
with $\myR \sHom_{\sO_{\Sigma}}(\DuBois{\Sigma}, \omega_{\Sigma}^{\mydot})$ and the
middle object can be identified with $\myR \sHom_{\sO_{X}}(\DuBois{X},
\omega_{X}^{\mydot})$.  Recall that $X$ has Du~Bois singularities if and only if the
natural map $\myR \sHom_{\sO_{X}}(\DuBois{X}, \omega_{X}^{\mydot}) \rightarrow
\omega_X^{\mydot}$ is an isomorphism.  Therefore, the object $\pi_* \omega_{\tld
  X}(E)$ is closely related to whether or not $X$ has Du~Bois singularities.  This
inspired the following result which we state (but do not prove).

\begin{theorem} \cite[Theorem 3.1]{KSS08}
  \label{ThmCMDuBoisCriterion}
  Suppose that $X$ is normal and Cohen-Macaulay.  Let $\pi : X' \rightarrow X$ be a
  log resolution, and denote the reduced exceptional divisor of $\pi$ by $G$.  Then
  $X$ has Du~Bois singularities if and only if $\pi_* \omega_{X'}(G) \simeq
  \omega_X$.
\end{theorem}
\begin{proof}
  We will not prove this.  The main idea is to show that

\hfill $\pi_* \omega_{X'}(G) \simeq \myH^{-\dim X} \left( \myR
    \sHom_{\sO_X}(\DuBois{X}, \omega_X^{\mydot}) \right).$ \hfill
\end{proof}

Related results can also be obtained in the non-normal Cohen-Macaulay case, see
\cite{KSS08} for details.

\begin{remark}
  The submodule $\pi_* \omega_{X'}(G) \subseteq \omega_X$ is independent of the
  choice of log resolution.  Thus this submodule may be viewed as an invariant which
  partially measures how far a scheme is from being Du~Bois (compare with
  \cite{FujinoNonLCSheaves}).
\end{remark}

As an easy corollary, we obtain another proof that rational singularities are Du~Bois
(this time via the Kempf-criterion for rational singularities).
\begin{corollary}
  If $X$ has rational singularities, then $X$ has Du~Bois singularities.
\end{corollary}
\begin{proof}
  Since $X$ has rational singularities, it is Cohen-Macaulay and normal.
  Then $\pi_* \omega_{X'} = \omega_X$ but we also have $\pi_* \omega_{X'} \subseteq
  \pi_* \omega_{X'}(G) \subseteq \omega_X$, and thus $\pi_* \omega_{X'}(G) =
  \omega_X$ as well.  Then use Theorem \ref{ThmCMDuBoisCriterion}.
\end{proof}

We also see immediately that log canonical singularities coincide with Du~Bois
singularities in the Gorenstein case.
\begin{corollary}
\label{CorGorLCImpliesDuBois}
Suppose that $X$ is Gorenstein and normal.  Then $X$ is Du~Bois if and only if $X$ is
log canonical.
\end{corollary}
\begin{proof}
  $X$ is easily seen to be log canonical if and only if $\pi_* \omega_{X'/X}(G) \simeq
  \sO_X$.  The projection formula then completes the proof.
\end{proof}

In fact, a slightly jazzed up version of this argument can be used to show that every
Cohen-Macaulay log canonical pair is Du~Bois, see \cite[Theorem 3.16]{KSS08}.

\section{The Koll\'ar-Kov\'acs splitting criterion}

Recently Koll\'ar and the first named author found a rather flexible criterion for
Du~Bois singularities.

\begin{thm}\cite{KK09a}
  \label{thm:db-criterion}
  Let ${f}: Y\to X$ be a proper 
  morphism between reduced schemes of finite type over $\bC$, $W\subseteq X$ an
  arbitrary subscheme, and $F\leteq {f}^{-1}(W)$, equipped with the induced reduced
  subscheme structure.  Let $\sI_{W\subseteq X}$ denote the ideal sheaf of $W$ in $X$
  and $\sI_{F\subseteq Y}$ the ideal sheaf of $F$ in $Y$.  Assume that the natural
  map $\varrho$
  $$
  \xymatrix{ \sI_{W\subseteq X} \ar[r]_-\varrho & \myR{f}_*\sI_{F\subseteq Y}
    \ar@{-->}@/_1.5pc/[l]_{\varrho'} }
  $$
  admits a left inverse $\varrho'$, that is, $\rho'\circ\rho=\id_{\sI_{W\subseteq
      X}}$. Then if $Y,F$, and $W$ all have DB singularities, then so does $X$.
\end{thm}

For the proof, please see the original paper.

\begin{subrem}
  Notice that it is not required that ${f}$ be birational. On the other hand the
  assumptions of the theorem and \cite[Thm~1]{KovacsRat} imply that if $Y\setminus F$
  has rational singularities, e.g., if $Y$ is smooth, then $X\setminus W$ has
  rational singularities as well.
\end{subrem}

This theorem is used to derive various consequences in \cite{KK09a}, some of which
are formally unrelated to Du~Bois singularities. We will mention some of these in the
sequel, but the interested reader should look at the original article to obtain the
full picture.

\section{Log canonical singularities are Du~Bois}

Log canonical and Du~Bois singularities are very closely related as we have seen in
the previous sections.  This was first observed in \cite{Ishii85}, see also
\cite{WatanabeKimioPlurigenera} and \cite{Ishii87b}.

Recently, Koll\'ar and the first named author gave a proof that log canonical
singularities are Du~Bois using Theorem~\ref{thm:db-criterion}.  We will sketch some
ideas of the proof here.  There are two main steps.  First, one shows that the
non-klt locus of a log canonical singularity is Du~Bois (this generalizes
\cite{AmbroSeminormalLocus} and \cite[Corollary 7.3]{SchwedeCentersOfFPurity}).  Then
one uses Theorem~\ref{thm:db-criterion} to show that this property is enough to
conclude that $X$ itself is Du~Bois.  For the first part we refer the reader to the
original paper. The key point of the second part is contained in the following Lemma.
Here we give a different proof than in \cite{KK09a}.

\begin{lemma}
  Suppose $(X, \Delta)$ is a log canonical pair and that the reduced non-klt locus of
  $(X, \Delta)$ has Du~Bois singularities.  Then $X$ has Du~Bois singularities.
\end{lemma}
\begin{proof}
  First recall that the multiplier ideal $\mJ(X, \Delta)$ is precisely the defining
  ideal of the non-klt locus of $(X, \Delta)$ and since $(X, \Delta)$ is log
  canonical, it is a radical ideal.  We set $\Sigma \subseteq X$ to be the reduced
  subscheme of $X$ defined by this ideal.  Since the statement is local, we may
  assume that $X$ is affine and thus that $X$ is embedded in a smooth scheme $Y$.  We
  let $\pi : \tld Y \rightarrow Y$ be an embedded resolution of $(X, \Delta)$ in $Y$
  and we assume that $\pi$ is an isomorphism outside the singular locus of $X$.  Set
  $\overline \Sigma$ to be the reduced-preimage of $\Sigma$ (which we may assume is a
  divisor in $\tld Y$) and let $\tld X$ denote the strict transform of $X$.  We
  consider the following diagram of exact triangles.
  \[
  \xymatrix{ & A^{\mydot} \ar[d]^{\alpha} \ar[r] & B^{\mydot} \ar[d]^{\beta} \ar[r] &
    C^{\mydot} \ar[d]^{\gamma} \ar[r]^-{+1} & \\
    0 \ar[r] & \mJ(X, \Delta) \ar[d]^{\delta} \ar[r] & \sO_X \ar[d] \ar[r] &
    \sO_{\Sigma} \ar[d]^{\varepsilon} \ar[r] & 0 \\
    & \myR \pi_* \sO_{\tld X}(-\overline \Sigma) \ar[r] & \myR \pi_* \sO_{\overline
      \Sigma \cup \tld X} \ar[r] & \myR \pi_* \sO_{\overline \Sigma} \ar[r]^-{+1} & \\
  }
  \]
  Here the first row is made up of objects in $D^{b}_{\coherent}(X)$ needed to make
  the columns into exact triangles.  Since $\Sigma$ has Du~Bois singularities, the
  map $\varepsilon$ is an isomorphism and so $C^{\mydot} \simeq 0$.  On the other
  hand, there is a natural map $\myR \pi_* \sO_{\tld X}(-\overline \Sigma)
  \rightarrow \myR \pi_* \sO_{\tld X}(K_{\tld X} - \pi^*(K_X + \Delta)) \simeq \mJ(X,
  \Delta)$ since $(X, \Delta)$ is log canonical.  This implies that the map $\alpha$
  is the zero map in the derived category.  However, we then see that $\beta$ is also
  zero in the derived category which implies that $\sO_X \rightarrow \myR \pi_*
  \sO_{\overline \Sigma \cup \tld X}$ has a left inverse.  Therefore, $X$ has Du~Bois
  singularities (since $\overline \Sigma \cup \tld X$ has simple normal crossing
  singularities) by Theorems~\ref{ThmDuBoisSplitting} and \ref{EasyDuBoisCriterion}.
\end{proof}

\section{Applications to moduli spaces and vanishing theorems}

The connection between log canonical and Du~Bois singularities have many useful
applications in moduli theory. We will list a few without proof.

\begin{demo}{Setup}\label{assume}%
  Let $\phi:X\to B$ be a flat projective morphism of complex varieties with $B$
  connected.  Assume that for all $b\in B$ there exists a $\bQ$-divisor $D_b$ on
  $X_b$ such that $(X_b,D_b)$ is log canonical.
\end{demo}

\begin{remark}
  Notice that it is not required that the divisors $D_b$ form a family.
\end{remark}

\begin{thm}\cite{KK09a}
  \label{thm:coh-inv}
  Under the assumptions in \eqref{assume}, $h^i(X_b, \sO_{X_b})$ is independent of
  $b\in B$ for all $i$.
\end{thm}

\begin{thm}\cite{KK09a}
  \label{thm:cm-inv}
  Under the assumptions in \eqref{assume} if one fiber of $\phi$ is Cohen-Macaulay
  (resp.\ $S_k$ for some $k$), then all fibers are Cohen-Macaulay (resp.\ $S_k$).
\end{thm}

\begin{thm}\cite{KK09a}\label{thm:flat}
  Under the assumptions in \eqref{assume} the cohomology sheaves
  $h^{i}(\omega_\phi^\mydot)$ are flat over $B$, where $\omega_\phi^\mydot$ denotes
  the {relative dualizing complex} of $\phi$.
\end{thm}

Du~Bois singularities also appear naturally in vanishing theorems.  As a culmination
of the work of Tankeev, Ramanujam, Miyaoka, Kawamata, Viehweg, Koll\'ar,
and Esnault-Viehweg, Koll\'ar proved a rather general form of a Kodaira-type vanishing
theorem in \cite[9.12]{KollarShafarevich}.  Using the same ideas this was slightly
generalized to the following theorem in \cite{KSS08}.

\begin{theorem}[{\cite[9.12]{KollarShafarevich}, \cite[6.2]{KSS08}}]
  \label{thm:several}
  Let $X$ be a proper variety and $\sL$ a line bundle on $X$. Let $\sL^m\simeq \sO_X
  (D)$, where $D=\sum d_i D_i$ is an effective divisor, and let $s$ be a global
  section whose zero divisor is $D$.  Assume that $0<d_i<m$ for every $i$. Let $Z$ be
  the scheme obtained by taking the $m^{th}$ root of $s$ (that is, $Z = X[\sqrt{s}]$
  using the notation from \cite[9.4]{KollarShafarevich}). Assume further that
  \begin{equation*}
    H^j(Z, \bC_Z) \to H^j(Z, \sO_Z)
  \end{equation*}
  is surjective. Then for any collection of $b_i\geq 0$ the natural map
  $$
  H^j\left(X, \sL^{-1}\left(-\sum b_iD_i\right)\right) \to H^j(X, \sL^{-1})
  $$
  is surjective.
\end{theorem}

This, combined with the fact that log canonical singularities are Du~Bois yields that
Kodaira vanishing holds for log canonical pairs:

\begin{thm}\cite[6.6]{KSS08}
  \label{KodairaVanishingForLC}
  Kodaira vanishing holds for Cohen-Macaulay  semi-log canonical varieties: Let
  $(X,\Delta)$ be a projective Cohen-Macaulay  semi-log canonical pair and
  $\sL$ an ample line bundle on $X$. Then $\coh i.X.\sL^{-1}.=0$ for $i<\dim X$.
\end{thm}

It turns out that Du~Bois singularities appear naturally in other kinds of vanishing
theorems. Let us cite one of those here.

\begin{thm}\cite[9.3]{GKKP09}
  \label{cor:idealsheafvanishingwithboundary}
  Let $(X, D)$ be a log canonical reduced pair of dimension $n \geq 2$, $\pi :
  \widetilde X \to X$ a log resolution with $\pi$-exceptional set $E$, and
  $\widetilde D = \Supp \bigl( E + \pi^{-1}D \bigr)$. Then
  $$ R^ {n-1}\pi_*\sO_{\widetilde X}(- \widetilde D) = 0.$$
\end{thm}

\section{Deformations of Du~Bois singularities}

Given the importance of Du~Bois singularities in moduli theory it is an important
obvious question whether they are invariant under small deformation.

It is relatively easy to see from the construction of the Du~Bois complex that a
general hyperplane section (or more generally, the general member of a base point
free linear system) on a variety with Du~Bois singularities again has Du~Bois
singularities. Therefore the question of deformation follows from the following.

\begin{conj}(cf.\ \cite{SteenbrinkMixed}) %
  \label{conj:DB-defo} %
  Let $D\subset X$ be a reduced Cartier divisor and assume that $D$
  has only Du~Bois singularities in a neighborhood of a point $x\in D$. Then $X$ has
  only Du~Bois singularities in a neighborhood of the point $x$.
\end{conj}

This conjecture was proved for isolated Gorenstein singularities by Ishii
\cite{Ishii86}.  Also note that rational singularities satisfy this property, see
\cite{ElkikDeformationsOfRational}.

We also have the following easy corollary of the results presented earlier:

\begin{thm}\label{thm:Gor-defo}
  Assume that $X$ is Gorenstein and $D$ is normal.\footnote{this condition is
    actually not necessary, but the proof becomes rather involved without it.}  Then
  the statement of Conjecture~\ref{conj:DB-defo} is true.
\end{thm}

\begin{proof}
  The question is local so we may restrict to a neighborhood of $x$. %
  If $X$ is Gorenstein, then so is $D$ as it is a Cartier divisor. Then $D$ is log
  canonical by \eqref{CorGorLCImpliesDuBois}, and then the pair $(X,D)$ is also log
  canonical by inversion of adjunction \cite{KawakitaInversion}. (Recall that if $D$
  is normal, then so is $X$ along $D$).  This implies that $X$ is also log canonical
  and thus Du~Bois.
\end{proof}

It is also stated in \cite[3.2]{KovacsDuBoisLC2} that the conjecture holds in full
generality. Unfortunately, the proof is not complete. The proof published there works
if one assumes that the non-Du~Bois locus of $X$ is contained in $D$. For instance,
one may assume that this is the case if the non-Du~Bois locus is isolated.

The problem with the proof is the following: it is stated that by taking hyperplane
sections one may assume that the non-Du~Bois locus is isolated. However, this is
incorrect. One may only assume that the \emph{intersection} of the non-Du~Bois locus
of $X$ with $D$ is isolated. If one takes a further general section then it will miss
the intersection point and then it is not possible to make any conclusions about that
case.

Therefore currently the best known result with regard to this conjecture is the
following:

\begin{thm}\cite[3.2]{KovacsDuBoisLC2}
  Let $D\subset X$ be a reduced Cartier divisor and assume that $D$ has only Du~Bois
  singularities in a neighborhood of a point $x\in D$ and that $X\setminus D$ has
  only Du~Bois singularities.  Then $X$ has only Du~Bois singularities in a
  neighborhood of $x$.
\end{thm}

Experience shows that divisors not in general position tend to have worse
singularities then the ambient space in which they reside. Therefore one would in
fact expect that if $X\setminus D$ is reasonably nice, and $D$ has Du~Bois
singularities, then perhaps $X$ has even better ones.

We have also seen that rational singularities are Du~Bois and at least Cohen-Macaulay
Du~Bois singularities are not so far from being rational cf.\
\ref{ThmCMDuBoisCriterion}. The following result of the second named author supports
this philosophical point.

\begin{thm}\cite[Thm.\ 5.1]{SchwedeEasyCharacterization}
  Let $X$ be a reduced scheme of finite type over a field of characteristic zero, $D$
  a Cartier divisor that has Du~Bois singularities and assume that $X\setminus D$ is
  smooth.  Then $X$ has rational singularities (in particular, it is Cohen-Macaulay).
\end{thm}

Let us conclude with a conjectural generalization of this statement:

\begin{conj}\label{conj:DB-to-rtl}
  Let $X$ be a reduced scheme of finite type over a field of characteristic zero, $D$
  a Cartier divisor that has Du~Bois singularities and assume that $X\setminus D$ has
  rational singularities.  Then $X$ has rational singularities (in particular, it is
  Cohen-Macaulay).
\end{conj}

Essentially the same proof as in \eqref{thm:Gor-defo} shows that this is also true
under the same additional hypotheses.

\begin{thm}
  Assume that $X$ is Gorenstein and $D$ is normal.\footnote{again, this condition is
    not necessary, but the proof becomes rather involved without it} Then the
  statement of Conjecture~\ref{conj:DB-to-rtl} is true.
\end{thm}

\begin{proof}
  If $X$ is Gorenstein, then so is $D$ as it is a Cartier divisor. Then by
  \eqref{CorGorLCImpliesDuBois} $D$ is log canonical. Then by inversion of adjunction
  \cite{KawakitaInversion} the pair $(X,D)$ is also log canonical near $D$. (Recall
  that if $D$ is normal, then so is $X$ along $D$).

  As $X$ is Gorenstein and $X\setminus D$ has rational singularities, it follows that
  $X \setminus D$ has canonical singularities.  Then $X$ has only canonical
  singularities everywhere. This can be seen by observing that $D$ is a Cartier
  divisor and examining the discrepancies that lie over $D$ for $(X, D)$ as well as
  for $X$. Therefore, by \cite{ElkikRationalityOfCanonicalSings}, $X$ has only
  rational singularities along $D$.
\end{proof}

\section{Characteristic $p > 0$ analogs of Du~Bois singularities}

Starting in the early 1980s, the connections between singularities defined by the
action of the Frobenius morphism in characteristic $p > 0$ and singularities defined
by resolutions of singularities started to be investigated, cf.\
\cite{FedderFPureRat}.  After the introduction of tight closure in
\cite{HochsterHunekeTC1}, a precise correspondence between several classes of
singularities was established.  See, for example, \cite{FedderWatanabe},
\cite{MehtaSrinivasFPureSurface}, \cite{HaraWatanabeFRegFPure},
\cite{SmithFRatImpliesRat}, \cite{HaraRatImpliesFRat},
\cite{MehtaSrinivasRatImpliesFRat}, \cite{SmithMultiplierTestIdeals},
\cite{HaraACharacteristicPAnalogOfMultiplierIdealsAndApplications},
\cite{HaraYoshidaGeneralizationOfTightClosure},
\cite{TakagiInterpretationOfMultiplierIdeals}, \cite{TakagiWatanabeFPureThresh},
\cite{TakagiPLTAdjoint}.  The second name author partially extended this
correspondence in his doctoral dissertation by linking Du~Bois singularities with
$F$-injective singularities, a class of singularities defined in
\cite{FedderFPureRat}.  The currently known implications are summarized below.
\vskip 1em
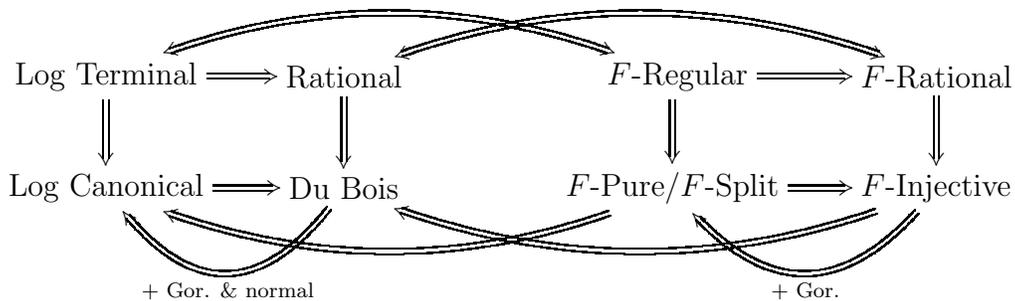
\[
\xymatrix{
  \text{Log Terminal} \ar@/^2pc/@{<=>}[rrr] \ar@{=>}[r] \ar@{=>}[d]& \text{Rational}
  \ar@{=>}[d] \ar@/^2pc/@{<=>}[rrr] & & \text{
    $F$-Regular} \ar@{=>}[r] \ar@{=>}[d] & \text{$F$-Rational} \ar@{=>}[d]\\
  \text{Log Canonical} \ar@2{=>}[r] \ar@/_2pc/@{<=}[rrr]
  \ar@/_2.7pc/@{<=}[r]_-{\text{+ Gor. \& normal}} & \text{Du~Bois}
  \ar@/_2pc/@{<=}[rrr] & & \text{$F$-Pure/$F$-Split}
  \ar@{=>}[r] \ar@/_2.7pc/@{<=}[r]_-{\text{+ Gor.}} & \text{$F$-Injective}\\
  & & & & \\
}
\]
We will give a short proof that normal Cohen-Macaulay singularities of dense
$F$-injective type are Du~Bois, based on the characterization of Du~Bois
singularities given in Section \ref{SectionAMultiplierIdealLikeCharacterization}.

Note that Du~Bois and $F$-injective singularities also share many common properties.
For example $F$-injective singularities are also seminormal \cite[Theorem
4.7]{SchwedeFInjectiveAreDuBois}.

First however, we will define $F$-injective singularities (as well as some necessary prerequisites).

\begin{definition}
  \label{DefinitionFFinite} Suppose that $X$ is a scheme of characteristic $p > 0$
  with absolute Frobenius map $F : X \rightarrow X$.  We say that $X$ is \emph{$F$-finite} if
  $F_* \sO_X$ is a coherent $\sO_X$-module.  A ring $R$ is called \emph{$F$-finite} if
  the associated scheme $\Spec R$ is $F$-finite.
\end{definition}

\begin{remark}
  Any scheme of finite type over a perfect field is $F$-finite, see for example
  \cite{FedderFPureRat}.
\end{remark}

\begin{definition}
  Suppose that $(R, \bm)$ is an $F$-finite local ring.  We say that $R$ is
  \emph{$F$-injective} if the induced Frobenius map $\xymatrix{H^i_{\bm}(R)
    \ar[r]^{F} & H^i_{\bm}(R)}$ is injective for every $i > 0$.  We say that an
  $F$-finite scheme is \emph{$F$-injective} if all of its stalks are $F$-injective
  local rings.
\end{definition}

\begin{remark}
  If $(R, \bm)$ is $F$-finite, $F$-injective and has a dualizing complex, then $R_{Q}$ is also $F$-injective for any $Q \in \Spec R$.  This
  follows from local duality, see \cite[Proposition 4.3]{SchwedeFInjectiveAreDuBois} for details.
\end{remark}

\begin{lemma}
  Suppose $X$ is a Cohen-Macaulay scheme of finite type over a perfect field $k$.
  Then $X$ is $F$-injective if and only if the natural map $F_* \omega_X \rightarrow
  \omega_X$ is surjective.
\end{lemma}
\begin{proof}
  Without loss of generality (since $X$ is Cohen-Macaulay) we can assume that $X$ is
  equidimensional.  Set $f : X \rightarrow \Spec k$ to be the structural morphism.
  Since $X$ is finite type over a perfect field, it has a dualizing complex
  $\omega_X^{\mydot} = f^! k$ and we set $\omega_X = h^{-\dim X}(\omega_X^{\mydot})$.
  Since $X$ is Cohen-Macaulay, $X$ is $F$-injective if and only if the Frobenius map
  $\xymatrix{H^{\dim X}_{x}(\sO_{X, x}) \ar[r] & H^{\dim X}_{x}(F_* \sO_{X, x})}$ is
  injective for every closed point $x \in X$.  By local duality, see
  \cite[Theorem 6.2]{HartshorneResidues} or \cite[Section 3.5]{BrunsHerzog}, such a map is injective if
  and only if the dual map $F_* \omega_{X, x} \rightarrow \omega_{X, x}$ is
  surjective.  But that map is surjective, if and only if the map of sheaves $F_*
  \omega_X \rightarrow \omega_X$ is surjective.
\end{proof}

We now briefly describe reduction to characteristic $p > 0$.  Excellent and far more
complete references include \cite[Section
2.1]{HochsterHunekeTightClosureInEqualCharactersticZero} and
\cite[II.5.10]{Kollar96}.  Also see \cite{SmithTightClosure} for a more elementary
introduction.

Let $R$ be a finitely generated algebra over a field $k$ of characteristic zero.
Write $R = k[x_1, \ldots, x_n]/I$ for some ideal $I$ and let $S$ denote $k[x_1,
\ldots, x_n]$.  Let $X = \Spec R$ and $\pi : \tld X \rightarrow X$ a log resolution
of $X$ corresponding to the blow-up of an ideal $J$.  Let $E$ denote the reduced
exceptional divisor of $\pi$. Then $E$ is the subscheme defined by the radical of the
ideal $J \cdot \O_{\tld X}$.

There exists a finitely generated $\bZ$-algebra $A \subset k$ that includes all the
coefficients of the generators of $I$ and $J$, a finitely generated $A$ algebra $R_A
\subset R$, an ideal $J_A \subset R_A$, and schemes $\tld X_A$ and $E_A$ of finite
type over $A$ such that $R_A \tensor_A k = R$, $J_A R = J$, $\tld X_A \times_{\Spec
  A} \Spec k = X$ and $E_A \times_{\Spec A} \Spec k = E$ with $E_A$ an effective
divisor with support defined by the ideal $J_A \cdot \O_{\tld X_A}$.  We may localize
$A$ at a single element so that $Y_A$ is smooth over $A$ and $E_A$ is a reduced
simple normal crossings divisor over $A$.  By further localizing $A$ (at a single
element), we may assume any finite set of finitely generated $R_A$ modules is
$A$-free (see \cite[3.4]{HunekeTightClosureBook} or
\cite[2.3]{HochsterRobertsFrobeniusLocalCohomology}) and we may assume that $A$
itself is regular.  We may also assume that a fixed affine cover of $E_A$ and a fixed
affine cover of $\tld X_A$ are also $A$-free.

We will now form a family of positive characteristic models of $X$ by looking at all
the rings $R_t = R_A \tensor_A k(t)$ where $k(t)$ is the residue field of a maximal
ideal $t \in T = \Spec A$.  Note that $k(t)$ is a finite, and thus perfect, field of
characteristic $p$.  We may also tensor the various schemes $X_A$, $E_A$, etc. with
$k(t)$ to produce a characteristic $p$ model of an entire situation.

By making various cokernels of maps free $A$-modules, we may also assume that maps
between modules that are surjective (respectively injective) over $k$ correspond to
surjective (respectively injective) maps over $A$, and thus surjective (respectively
injective) in our characteristic $p$ model as well; see
\cite{HochsterHunekeTightClosureInEqualCharactersticZero} for details.

\begin{definition}
  A ring $R$ of characteristic zero is said to have dense $F$-injective type if for
  every family of characteristic $p \gg 0$ models with $A$ chosen sufficiently large,
  we have that a Zariski dense set of those models (over $\Spec A$) have
  $F$-injective singularities.
\end{definition}

\begin{theorem}\cite{SchwedeFInjectiveAreDuBois}
  \label{thm:F-injctive-is-DB}
  Let $X$ be a reduced scheme of finite type over $\bC$ and assume that it has dense
  $F$-injective type.  Then $X$ has Du~Bois singularities.
\end{theorem}
\begin{proof}
  We only provide a proof in the case that $X$ is normal and Cohen-Macaulay.  For a
  complete proof, see \cite{SchwedeFInjectiveAreDuBois}.  Let $\pi : \tld X
  \rightarrow X$ be a log resolution of $X$ with exceptional divisor $E$.  We reduce
  this entire setup to characteristic $p \gg 0$ such that the corresponding $X$ is
  $F$-injective.  Let $F^e : X \rightarrow X$ be the $e$-iterated Frobenius map.

  We have the following commutative diagram,
  \[
  \xymatrix{ F^e_* \pi_* \omega_{\tld{X}}(p^e E) \ar[d]_{\rho} \ar[r] & \pi_*
    \omega_{\tld{X}}(E) \ar[d]^{\beta} \\
    F^e_* \omega_X \ar[r]^{\phi} & \omega_X \\
  }
  \]
  where the horizontal arrows are induced by the dual of the Frobenius map, $\sO_X
  \rightarrow F^e_* \sO_X$ and the vertical arrows are the natural maps induced by
  $\pi$.  By hypothesis, $\phi$ is surjective.  On the other hand, for $e > 0$
  sufficiently large, the map labeled $\rho$ is an isomorphism.  Therefore the map
  $\phi \circ \rho$ is surjective which implies that the map $\beta$ is also
  surjective.  But as this holds for a dense set of primes, it must be surjective in
  characteristic zero as well, and in particular, as a consequence $X$ has Du~Bois
  singularities.
\end{proof}

It is not known whether the converse of this statement is true:

\begin{demo}{\bf Open Problem}
  If $X$ has Du~Bois singularities, does it have dense F-injective type?
\end{demo}

Since $F$-injective singularities are known to be closely related to Du~Bois
singularities, it is also natural to ask how $F$-injective singularities deform cf.\
Conjecture~\ref{conj:DB-defo}.  In general, this problem is also open.
\begin{demo}{\bf Open Problem}
  If a Cartier divisor $D$ in $X$ has $F$-injective singularities, does $X$ have
  $F$-injective singularities near $D$?
\end{demo}
In the case that $X$ (equivalently $D$) is Cohen-Macaulay, the answer is affirmative,
see \cite{FedderFPureRat}.  In fact, Fedder defined $F$-injective singularities
partly because they seemed to deform better than $F$-pure singularities (the
conjectured analog of log canonical singularities).

\appendix{Connections with Buchsbaum rings}

In this section we discuss the links between Du~Bois singularities and Buchsbaum
rings.  Du~Bois singularities are not necessarily Cohen-Macaulay, but in many cases,
they are Buchsbaum (a weakening of Cohen-Macaulay).

Recall that a local ring $(R, \bm, k)$ has \emph{quasi-Buchsbaum} singularities if
$\bm H^i_{\bm}(R) = 0$ for all $i < \dim R$.  Further recall that a ring is called
\emph{Buchsbaum} if $\tau^{\dim R} \myR \Gamma_{\bm}(R)$ is quasi-isomorphic to a
complex of $k$-vector spaces.  Here $\tau^{\dim R}$ is the brutal truncation of the complex at the
$\dim R$ location.  Note that this is not the usual definition of Buchsbaum
singularities, rather it is the so-called Schenzel's criterion, see
\cite{SchenzelApplicationsOfDualizingComplexes}.  Notice that Cohen-Macaulay
singularities are Buchsbaum (after truncation, one obtains the zero-object in the derived category).

It was proved by Tomari that isolated Du~Bois singularities are quasi-Buchsbaum (a
proof can be found in \cite[Proposition 1.9]{Ishii85}), and then by Ishida that
isolated Du~Bois singularities were in fact Buchsbaum.  Here we will briefly review
the argument to show that isolated Du~Bois singularities are quasi-Buchsbaum since
this statement is substantially easier.

\begin{proposition}
  Suppose that $(X, x)$ is an isolated Du~Bois singularity with $R = \sO_{X, x}$.
  Then $R$ is quasi-Buchsbaum.
\end{proposition}
\begin{proof}
  Note that we may assume that $X$ is affine.  Since $\Spec R$ is regular
  outside its the maximal ideal $\bm$, it is clear that some power of $\bm$
  annihilates $H^i_{\bm}(R)$ for all $i < \dim R$.  We need to show that the smallest
  power for which this happens is $1$.  We let $\pi : \tld X \rightarrow X$ be a log
  resolution with exceptional divisor $E$ as in Theorem
  \ref{TheoremSteenbrinkIsolatedCharacterization}.  Since $X$ is affine, we see that
  $H^i_{\bm}(R) \simeq H^{i-1}(X \setminus \{\bm\}, \sO_X) \simeq H^{i-1}(\tld X
  \setminus E, \sO_{\tld X})$ for all $i > 0$. Therefore, it is enough to show that
  $\bm H^{i-1}(\tld X \setminus E, \sO_{\tld X}) = 0$ for all $i < \dim X$.  In other
  words, we need to show that $\bm H^{i}(\tld X \setminus E, \sO_{\tld X}) = 0$ for
  all $i < \dim X - 1$.

  We examine the long exact sequence
  \[
  \xymatrix@C=15pt{ \ldots \ar[r] & H^{i-1}(\tld X \setminus E, \sO_{\tld X}) \ar[r]
    & H^i_{E}(\tld X, \sO_{\tld X}) \ar[r] & H^i(\tld X, \sO_{\tld X}) \ar[r] &
    H^{i}(\tld X \setminus E,
    \sO_{\tld X}) \ar[r] & \ldots \\
  }
  \]
  Now, $H^i_{E}(\tld X, \sO_{\tld X}) = \myR^i( \Gamma_{\bm} \circ \pi_* )\left(
    \sO_X \right)$ which vanishes for $i < \dim X$ by the Matlis dual of
  Grauert-Riemenschneider vanishing.  Therefore $H^{i}(\tld X \setminus E, \sO_{\tld
    X}) \simeq H^i(\tld X, \sO_{\tld X})$ for $i < \dim X - 1$.  Finally, since $X$
  is Du~Bois, $H^i(\tld X, \sO_{\tld X}) = H^i(E, \sO_E)$ by Theorem
  \ref{TheoremSteenbrinkIsolatedCharacterization}.  But it is obvious that $\bm
  H^i(E, \sO_E) = 0$ since $E$ is a \emph{reduced} divisor whose image in $X$ is the
  point corresponding to $\bm$.  The result then follows.
\end{proof}

It is also easy to see that isolated $F$-injective singularities are also quasi-Buchsbaum.

\begin{proposition}
  Suppose that $(R, \bm)$ is a local ring that is $F$-injective.  Further suppose
  that $\Spec R \setminus \{ \bm \}$ is Cohen-Macaulay.  Then $(R, \bm)$ is
  quasi-Buchsbaum.
\end{proposition}
\begin{proof}
  Since the punctured spectrum of $R$ is Cohen-Macaulay, $H^i_{\bm}(R)$ is
  annihilated by some power of $\bm$ for $i < \dim R$.  We will show that the
  smallest such power is $1$.  Choose $c \in \bm$.  Since $R$ is $F$-injective, $F^e
  : H^i_{\bm}(R) \rightarrow H^i_{\bm}(R)$ is injective for all $e > 0$.  Choose $e$
  large enough so that $c^{p^e} H^i_{\bm}(R)$ is zero for all $i < e$.  However, for
  any element $z \in H^i_{\bm}(R)$, $F^e(cz) = c^{p^e} F^e(z) \in c^{p^e}
  H^i_{\bm}(R) = 0$ for $i < \dim R$.  This implies that $cz = 0$ and so $\bm
  H^i_{\bm}(R) = 0$ for $i < \dim R$.
\end{proof}

Perhaps the most interesting open question in this area is the following:

\begin{demo}{\bf Open Problem}[(Takagi)]\label{open:F}
 Are $F$-injective singularities with isolated non-CM locus Buchsbaum?
\end{demo}

Given the close connection between $F$-injective and Du~Bois singularities,
this question naturally leads to the next one:

\begin{demo}{\bf Open Problem}\label{open:DB}
 Are Du~Bois singularities with isolated non-CM locus Buchsbaum?
\end{demo}

\appendix{Cubical hyperresolutions}
\label{AppendixOnHyper}

For the convenience of the reader we include a short appendix explaining the
construction of cubical hyperresolutions, as well as several examples.  We follow
\cite{GNPP88} and mostly use their notation.

First let us fix a small universe to work in.  Let $\uSch$ denote the category of
reduced schemes.  One should note that the usual fibred product of schemes $X
\times_S Y$ need not be reduced, even when $X$ and $Y$ are reduced.  We wish to
construct the fibred product in the category of reduced schemes.  Given any scheme
$W$ (reduced or not) with maps to $X$ and $Y$ over $S$, there is always a unique
morphism $W \rightarrow X \times_S Y$, which induces a natural unique morphism
$W_{\red} \rightarrow (X \times_S Y)_{\red}$.  It is easy to see that $(X \times_S
Y)_{\red}$ is the fibred product in the category of reduced schemes.

Let us denote by $\uOne$ the category $\{0\}$ and by $\uTwo$ the category $\{0
\rightarrow 1\}$.  Let $n$ be an integer $\geq -1$.  We denote by $\square^+_n$ the
product of $n+1$ copies of the category $\uTwo = \{0 \rightarrow 1\}$ \cite[I,
1.15]{GNPP88}.  The objects of $\square^+_n$ are identified with the sequences
$\alpha = (\alpha_0, \alpha_1, \ldots, \alpha_n)$ such that $\alpha_i \in \{0, 1\}$
for $0 \leq i \leq n$.  For $n = -1$, we set $\square^+_{-1} = \{ 0 \}$ and for $n =
0$ we have $\square^+_0 = \{0 \rightarrow 1\}$.  We denote by $\square_n$ the full
subcategory consisting of all objects of $\square^+_n$ except the initial object $(0,
\ldots, 0)$.  Clearly, the category $\square^+_n$ can be identified with the category
of $\square_n$ with an augmentation map to $\{0\}$.

\begin{definition}
  A \emph{diagram of schemes} is a functor $\Phi$ from a category $\sfC^{\opp}$ to
  the category of schemes.  A \emph{finite diagram of schemes} is a diagram of
  schemes such that the aforementioned category $\sfC$ has finitely many objects and
  morphisms; in this case such a functor will be called a \emph{$\sfC$-scheme}. A
  morphism of diagrams of schemes $\Phi : \sfC^{\opp} \rightarrow \uSch$ to $\Psi :
  \sfD^{\opp} \rightarrow \uSch$ is the combined data of a functor $\Gamma :
  \sfC^{\opp} \rightarrow \sfD^{\opp}$ together with a natural transformation of
  functors $\eta:\Phi \to \Psi \circ \Gamma$.
\end{definition}

\begin{remark}
  With the above definitions, the class of (finite) diagrams of schemes can be made
  into a category.  Likewise the set of $\sfC$-schemes can also be made into a
  category (where the functor $\Gamma : \sfC^{\opp} \rightarrow \sfC^{\opp}$ is
  always chosen to be the identity functor).
\end{remark}

\begin{remark}
  Let $I$ be a category.  If instead of a functor to the category of reduced schemes,
  one considers a functor to the category of topological spaces, or the category of
  categories, one can define $I$-topological spaces, and $I$-categories in the
  obvious way.
\end{remark}

If $X_{\mydot} : I^{\opp} \rightarrow \uSch$ is an $I$-scheme, and $i \in \Ob I$,
then $X_i$ will denote the scheme corresponding to $i$.  Likewise if $\phi \in \Mor
I$ is a morphism $\phi : j \rightarrow i$, then $X_\phi$ will denote the
corresponding morphism $X_\phi : X_i \rightarrow X_j$.  If $f : Y_{\mydot}
\rightarrow X_{\mydot}$ is a morphism of $I$-schemes, we denote by $f_i$ the induced
morphism $Y_i \rightarrow X_i$.  If $X_{\mydot}$ is an $I$-scheme, a closed
sub-$I$-scheme is a morphism of $I$-schemes $g : Z_{\mydot} \rightarrow X_{\mydot}$
such that for each $i \in I$, the map $g_i : Z_i \rightarrow X_i$ is a closed
immersion.  We will often suppress the $g$ of the notation if no confusion is likely
to arise.  More generally, any property of a morphism of schemes (projective, proper,
separated, closed immersion, etc...) can be generalized to the notion of a morphism
of $I$-schemes by requiring that for each object $i$ of $I$, $g_i$ has the desired
property (projective, proper, separated, closed immersion, etc...)

\begin{definition}\cite[I, 2.2]{GNPP88}
  Given a morphism of $I$-schemes $f : Y_{\mydot} \rightarrow X_{\mydot}$, we define
  the \emph{discriminant} \emph{of $f$} to be the smallest closed sub-$I$-scheme
  $Z_{\mydot}$ of $X_{\mydot}$ such that $f_i : (Y_i - (f_i^{-1}(Z_i)) ) \rightarrow
  (X_i - Z_i)$ is an isomorphism for all $i$.
\end{definition}

\begin{definition}\cite[I, 2.5]{GNPP88}
  Let $S_{\mydot}$ be an $I$-scheme, $f : X_{\mydot} \rightarrow S_{\mydot}$ a proper
  morphism of $I$-schemes, and $D_{\mydot}$ the discriminant of $f$.  We say that $f$
  is a \emph{resolution\footnote{A resolution is a \emph{distinct} notion from a
      cubic hyperresolution.} of $S_{\mydot}$} if $X_{\mydot}$ is a smooth $I$-scheme
  (meaning that each $X_i$ is smooth) and $\dim f_i^{-1} (D_i) < \dim S_i$, for all
  $i \in \ob I$.
\end{definition}

\begin{remark}
  This is the definition found in \cite{GNPP88}.  Note that the maps are not required
  to be surjective (of course, the ones one constructs in practice are usually
  surjective).

  Consider the following example: the map $k[x,y]/(xy) \rightarrow k[x]$ which sends
  $y$ to $0$.  I claim that the associated map of schemes is a ``resolution'' of the
  $*$-scheme, $\Spec k[x,y]/(xy)$.  The discriminant is $\Spec k[x,y]/(x)$.  The
  pre-image however is simply the origin on $k[x]$, which has lower dimension than
  ``1''.  Resolutions like this one are sometimes convenient to consider.

  On the other hand, this definition seems to allow something it perhaps shouldn't.
  Choose any variety $X$ of dimension greater than zero and a closed point $z \in X$.
  Consider the map $z \rightarrow X$ and consider the $*$-scheme $X$.  The
  discriminant is all of $X$.  However, the pre-image of $X$ is still just a point,
  which has lower dimension than $X$ itself, by hypothesis.

  In view of these remarks, sometimes it is convenient to assume also that $\dim D_i
  < \dim S_i$ for each $i \in \ob I$.  In the resolutions of $I$-schemes that we
  construct (in particular, in the ones that are used to that prove cubic
  hyperresolutions exist), this always happens.
\end{remark}

Let $I$ be a category.  The set of objects of $I$ can be given the following
pre-order relation, $i \leq j$ if and only if $\Hom_I(i,j)$ is nonempty.  We will say
that a category $I$ is ordered if this pre-order is a partial order and, for each $i
\in \ob I$, the only endomorphism of $i$ is the identity \cite[I, C, 1.9]{GNPP88}.
Note that a category $I$ is ordered if and only if all isomorphisms and endomorphisms
of $I$ are the identity.

It turns out of that resolutions of $I$-schemes always exist under reasonable
hypotheses.

\begin{theorem}  \cite[I, Theorem 2.6]{GNPP88}
\label{thmResolutionOfIScheme}
Let $S$ be an $I$-scheme of finite type over a field $k$.  Suppose that $k$ is a
field of characteristic zero and that $I$ is a finite ordered category. Then there
exists a resolution of $S$.
\end{theorem}

In order to construct a resolution $Y_{\mydot}$ of an $I$-scheme $X_{\mydot}$, it
might be tempting to simply resolve each $X_i$, set $Y_i$ equal to that resolution,
and somehow combine this data together.  Unfortunately this cannot work, as shown by
the example below.
\begin{example}
\label{PinchPointCannotBeResolved}
Consider the pinch point singularity,
\[
X = \Spec k[x,y,z]/(x^2 y - z^2) = \Spec k[s, t^2, st]
\]
and let $Z$ be the closed subscheme defined by the ideal $(s, st)$ (this is the
singular set).  Let $I$ be the category $\{0 \rightarrow 1\}$.  Consider the
$I$-scheme defined by $X_0 = X$ and $X_1 = Z$ (with the closed immersion as the map).
$X_1$ is already smooth, and if one resolves $X_0$, (that is, normalizes it) there is
no compatible way to map $X_1$ (or even another birational model of $X_1$) to it,
since its pre-image by normalization will be two-to-one onto $Z \subset X$!  The way
this problem is resolved is by creating additional components.  So to construct a
resolution $Y_{\mydot}$ we set $Y_1 = Z = X_1$ (since it was already smooth) and set
$Y_0 = \overline X_0 \coprod Z$ where $\overline X_0$ is the normalization of $X_0$.
The map $Y_1 \rightarrow Y_0$ just sends $Y_1$ (isomorphically) to the new component
and the map $Y_0 \rightarrow X_0$ is the disjoint union of the normalization and
inclusion maps.
\end{example}

One should note that although the theorem proving the existence of resolutions of
$I$-schemes is constructive, \cite{GNPP88}, it is often easier in practice to
construct an ad-hoc resolution.

Now that we have resolutions of $I$-schemes, we can discuss cubic hyperresolutions of
schemes, in fact, even diagrams of schemes have cubic hyperresolutions!  First we
will discuss a single iterative step in the process of constructing cubic
hyperresolutions.  This step is called a $2$-resolution.

\begin{definition} \cite[I, 2.7]{GNPP88}
\label{dfnTwoResolution}
Let $S$ be an $I$-scheme and $Z_{\mydot}$ a $\square^+_1 \times I$-scheme.  We say
that $Z_{\mydot}$ is a \emph{$2$-resolution} of $S$ if $Z_{\mydot}$ is defined by the
Cartesian square (pullback, or fibred product in the category of (reduced)
$I$-schemes) of morphisms of $I$-schemes below
\[
\xymatrix{
  Z_{11} \ar@{^{(}->}[r] \ar[d] & Z_{01} \ar[d]^f \\
  Z_{10} \ar@{^{(}->}[r] & Z_{00} \\
}
\]
where
\begin{itemize}
\item[i)]  $Z_{00} = S$.
\item[ii)]  $Z_{01}$ is a smooth $I$-scheme.
\item[iii)]  The horizontal arrows are closed immersions of $I$-schemes.
\item[iv)]  $f$ is a proper $I$-morphism
\item[v)] $Z_{10}$ contains the discriminant of $f$; in other words, $f$ induces an
  isomorphism of $(Z_{01})_i - (Z_{11})_i$ over $(Z_{00})_i - (Z_{10})_i$, for all $i
  \in \ob I$.
\end{itemize}
\end{definition}

Clearly $2$-resolutions always exist under the same hypotheses that resolutions of
$I$-schemes exist: set $Z_{01}$ to be a resolution, $Z_{10}$ to be discriminant (or
any appropriate proper closed sub-$I$-scheme that contains it), and $Z_{11}$ its
(reduced) pre-image in $Z_{01}$.

Consider the following example,

\begin{example}
  Let $I = \{0\}$ and let $S$ be the $I$-scheme $\Spec k[t^2, t^3]$.  Let $Z_{01} =
  \bA^1 = \Spec k[t]$ and $Z_{01} \rightarrow S = Z_{00}$ be the map defined by
  $k[t^2, t^3] \rightarrow k[t]$.  The discriminant of that map is the closed
  subscheme of $S = Z_{00}$ defined by the map $\phi: k[t^2, t^3] \rightarrow k$ that
  sends $t^2$ and $t^3$ to zero.  Finally we need to define $Z_{11}$.  The usual
  fibered product in the category of schemes is $k[t]/(t^2)$, but we work in the
  category of reduced schemes, so instead the fibered product is simply the
  associated reduced scheme (in this case $\Spec k[t]/(t)$).  Thus our $2$-resolution
  is defined by the diagram of rings pictured below.
\[
\xymatrix{
& k[t]/(t) & \\
k[t]/(t) \ar[ur] & & k[t] \ar[ul] \\
& \ar[ur] \ar[ul] k[t^2, t^3] & \\
}
\].
\end{example}

We need one more definition before defining a cubic hyperresolution,

\renewcommand\rd{\red}

\begin{definition} \cite[I, 2.11]{GNPP88} Let $r$ be an integer greater than or equal
  to $1$, and let $X^n_{\mydot}$ be a $\square^+_n \times I$-scheme, for $1 \leq n
  \leq r$.  Suppose that for all $n$, $1 \leq n \leq r$, the $\square^+_{n-1} \times
  I $-schemes $X^{n+1}_{00 \mydot}$ and $X^n_{1 \mydot}$ are equal.  Then we define,
  by induction on $r$, a $\square^+_r \times I$-scheme
  \[
  Z_{\mydot} = \rd(X_{\mydot}^1, X_{\mydot}^2, \ldots, X_{\mydot}^r)
  \]
  that we call the \emph{reduction} of $(X_{\mydot}^1, \ldots, X^r_{\mydot})$, in the
  following way: If $r = 1$, one defines $Z_{\mydot} = X_{\mydot}^1$, if $r = 2$ one
  defines $Z_{\mydot \mydot} = \rd(X_{\mydot}^1, X_{\mydot}^2)$ by
  \[
  Z_{\alpha \beta} = \Biggl\{
  \begin{array}{lcl}
    X^1_{0 \beta} & , &  \text{if } \alpha = (0, 0), \\
    X^2_{\alpha \beta} & , & \text{if } \alpha \in \square_1 , \\
  \end{array}
  \]
  for all $\beta \in \square^+_0$, with the obvious morphisms.  If $r > 2$, one
  defines $Z_{\mydot}$ recursively as $\rd(\rd(X_{\mydot}^1, \ldots,
  X_{\mydot}^{r-1}), X_{\mydot}^r)$.
\end{definition}

Finally we are ready to define cubic hyperresolutions.

\begin{definition} \cite[I, 2.12]{GNPP88} Let $S$ be an $I$-scheme.  A \emph{cubic
    hyperresolution augmented over $S$} is a $\square_r^+ \times I$-scheme
  $Z_{\mydot}$ such that
  \[
  Z_{\mydot} = \rd(X_{\mydot}^1, \ldots, X_{\mydot}^r),
  \]
  where

  \begin{enumerate}
  \item
    $X_{\mydot}^1$ is a $2$-resolution of $S$,
  \item
    for $1 \leq n < r$, $X_{\mydot}^{n+1}$ is a 2-resolution of $X_1^n$, and
  \item
    $Z_\alpha$ is smooth for all $\alpha \in \square_r$.
  \end{enumerate}
\end{definition}

Now that we have defined cubic hyperresolutions, we should note that they exist under reasonable hypotheses

\begin{theorem} \cite[I, 2.15]{GNPP88}
  \label{thmCubicHyperresolutionsExist}
  Let $S$ be an $I$-scheme.  Suppose that $k$ is a field of characteristic zero and
  that $I$ is a finite (bounded) ordered category.  Then there exists $Z_{\mydot}$, a
  cubic hyperresolution augmented over $S$ such that
  \[
  \dim Z_\alpha \leq \dim S - |\alpha| + 1, \forall \alpha \in \square_r .
  \]
\end{theorem}

Below are some examples of cubic hyperresolutions.

\begin{example}
  \label{CubicHyperresolutionsOfCurves}
  Let us begin by computing cubic hyperresolutions of curves so let $C$ be a curve.
  We begin by taking a resolution $\pi : \overline C \rightarrow C$ (where $\overline
  C$ is just the normalization).  Let $P$ be the set of singular points of $C$; thus
  $P$ is the discriminant of $\pi$.  Finally we let $E$ be the reduced exceptional
  set of $\pi$, therefore we have the following Cartesian square
  \[
  \xymatrix{
    E \ar[r] \ar[d] & \overline C \ar[d]^\pi \\
    P \ar[r] & C }
  \]
  It is clearly already a $2$-resolution of $C$ and thus a cubic-hyperresolution of
  $C$.
\end{example}

\begin{example}
  \label{CubicHyperresolutionsOfIsolatedPointsWithSmoothExceptionalSet}
  Let us now compute a cubic hyperresolution of a scheme $X$ whose singular locus is
  itself a smooth scheme, and whose reduced exceptional set of a strong resolution
  $\pi : \tld X \rightarrow X$ is smooth (for example, any cone over a smooth
  variety).  As in the previous example, let $\Sigma$ be the singular locus of $X$
  and $E$ the reduced exceptional set of $\pi$, Then the Cartesian square of reduced
  schemes
  \[
  \xymatrix{
    E \ar[r] \ar[d] & \tld X \ar[d]^\pi \\
    \Sigma \ar[r] & X \\
  }
  \]
  is in fact a $2$-resolution of $X$, just as in the case of curves above.
\end{example}

The obvious algorithm used to construct cubic hyperresolutions does not construct
hyperresolutions in the most efficient or convenient way possible.  For example,
applying the obvious algorithm to the intersection of three coordinate planes gives
us the following.

\begin{example}
  \label{exampleAlgorithmicThreeCoordinatePlanes}
  Let $X \cup Y \cup Z$ be the three coordinate planes in $\bA^3$.  In this example
  we construct a cubic hyperresolution using the obvious algorithm. What makes this
  construction different, is that the dimension is forced to drop when forming the
  discriminant of a resolution of a diagram of schemes.

  Yet again we begin the algorithm by taking a resolution and the obvious one is $\pi
  : (X \coprod Y \coprod Z) \rightarrow (X \cup Y \cup Z)$.  The discriminant is $B =
  (X \cap Y) \cup (X \cap Z) \cup (Y \cap Z)$, the three coordinate axes.  The fiber
  product making the square below Cartesian is simply the exceptional set $E = ((X
  \cap Y) \cup (X \cap Z)) \coprod ((Y \cap X) \cup (Y \cap Z)) \coprod ((Z \cap X)
  \cup (Z \cap Y))$.
  $$
  \xymatrix
  { %
    \scriptstyle E = \big((X \cap Y) \cup (X \cap Z)\big) \coprod \big((Y \cap X)
    \cup (Y \cap Z)\big) \coprod \big((Z \cap X) \cup (Z \cap Y)\big) \ar[r] \ar[d]^\phi
    & \scriptstyle \big(X \coprod Y \coprod Z\big) \ar[d]^\pi
    \\
    \scriptstyle {B = (X \cap Y) \cup (X \cap Z) \cup (Y \cap Z) \ar[r]} &
    \scriptstyle
    X \cup Y \cup Z\\
  }
  $$
  We now need to take a $2$-resolution of the $\uTwo$-scheme $\phi : E \rightarrow
  B$.  We take the obvious resolution that simply separates irreducible components.
  This gives us $\tld E \rightarrow \tld B$ mapping to $\phi : E \rightarrow B$.  The
  discriminant of $\tld E \rightarrow E$ is a set of three points $X_0$, $Y_0$ and
  $Z_0$ corresponding to the origins in $X$, $Y$ and $Z$ respectively.  The
  discriminant of the map $\tld B \rightarrow B$ is simply identified as the origin
  $A_0$ of our initial scheme $X \cup Y \cup Z$ (recall $B$ is the union of the three
  axes).  The union of that with the images of $X_0$, $Y_0$ and $Z_0$ is again just
  $A_0$.  The fiber product of the diagram
  \[
  \xymatrix{
    & (\tld E \rightarrow \tld B) \ar[d] \\
    (\{X_0, Y_0, Z_0 \} \ar[r] \rightarrow \{A_0\}) & (\phi : E \rightarrow B) \\
  }
  \]
  can be viewed as $\{Q_1, \ldots, Q_6\} \rightarrow \{P_1, P_2, P_3\}$ where $Q_1$
  and $Q_2$ are mapped to $P_1$ and so on (remember $E$ was the disjoint union of the
  coordinate axes of $X$, of $Y$, and of respectively $Z$, so $\tld E$ has six
  components and thus six origins).  Thus we have the following diagram
  
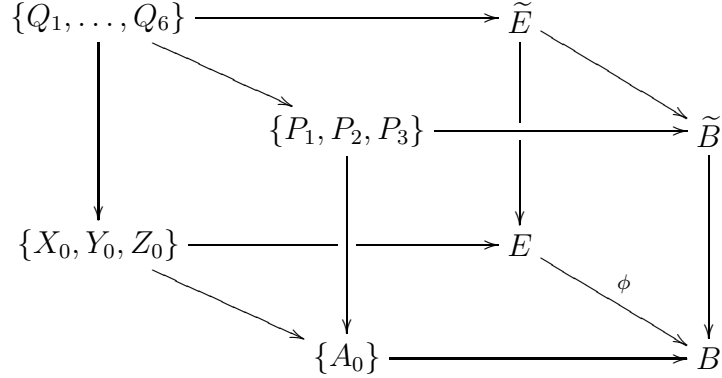
\[
  \xymatrix{%
    \{Q_1, \ldots, Q_6\} \ar@{->}[rr]\ar@{->}[dd] & & \tld E \ar@{->}'[d][dd]
    \\
    & \{P_1, P_2, P_3\} \ar@{<-}[ul]\ar@{->}[rrr]\ar@{->}[dd] & & & \tld B
    \ar@{<-}[ull]\ar@{->}[dd]
    \\
    \{X_0, Y_0, Z_0\} \ar@{->}'[r][rr] & & E
    \\
    & \{A_0\} \ar@{->}[rrr]\ar@{<-}[ul] & & & B \ar@{<-}[ull]_\phi
  }
  \]
  which we can combine with previous diagrams to construct a cubic hyperresolution.
\end{example}

\begin{remark}
  \label{exampleThreeCoordinatePlanes}
  It is possible to find a cubic hyperresolution for the three coordinate planes in
  $\bA^3$ in a different way.  Suppose that $S$ is the union of the three coordinate
  planes ($X$, $Y$, and $Z$) of $\bA^3$.  Consider the $\square_2$ or $\square^+_2$
  scheme defined by the diagram below (where the dotted arrows are those in
  $\square_2^+$ but not in $\square_2$).
  \[
  \qquad
   
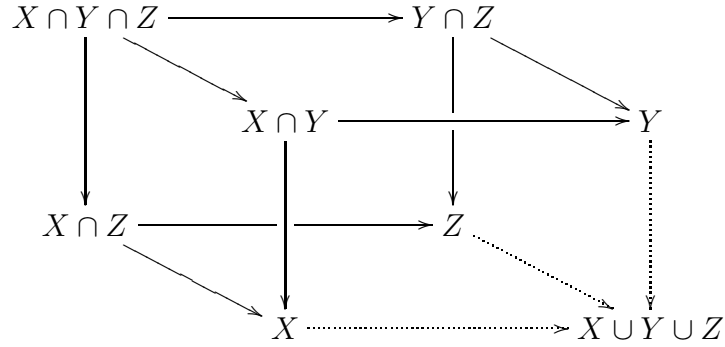
\xymatrix{ X \cap Y \cap Z \ar@{->}[rr]\ar@{->}[dd] & & Y \cap Z \ar@{->}'[d][dd]
    \\
    & X \cap Y \ar@{<-}[ul]\ar@{->}[rr]\ar@{->}[dd] & & Y \ar@{<-}[ul]\ar@{.>}[dd]
    \\
    X \cap Z \ar@{->}'[r][rr] & & Z
    \\
    & X \ar@{.>}[rr]\ar@{<-}[ul] & & X \cup Y \cup Z \ar@{<.}[ul]}
  \]
  One can verify that this is also a cubic hyperresolution of $X \cup Y \cup Z$.
\end{remark}

Now let us discuss sheaves on diagrams of schemes, as well as the related notions of push forward and its right derived functors.

\begin{definition} \cite[I, 5.3-5.4]{GNPP88} Let $X_{\mydot}$ be an $I$-scheme (or
  even an $I$-topological space).  We define a \emph{sheaf (or pre-sheaf) of abelian
    groups $F^{\mydot}$ on $X_{\mydot}$} to be the following data:
  \begin{enumerate}
  \item
    A sheaf (pre-sheaf) $F^i$ of abelian groups over $X_i$, for all $i \in \ob I$,
    and
  \item
    An $X_\phi$-morphism of sheaves $F^\phi : F^i \rightarrow (X_\phi)_* F^j$ for all
    morphisms $\phi : i \rightarrow j$ of $I$, required to be compatible in the
    obvious way.
  \end{enumerate}
\end{definition}

Given a morphism of diagrams of schemes $f_{\mydot} : X_{\mydot} \rightarrow
Y_{\mydot}$ one can construct a push-forward functor for sheaves on $X_{\mydot}$.

\begin{definition}\cite[I, 5.5]{GNPP88}
  Let $X_{\mydot}$ be an $I$-scheme, $Y_{\mydot}$ a $J$-scheme, $F^{\mydot}$ a sheaf
  on $X_{\mydot}$, and $f_{\mydot} : X_{\mydot} \rightarrow Y_{\mydot}$ a morphism of
  diagrams of schemes.  We define $(f_{\mydot})_* F^{\mydot}$ in the following way.
  For each $j \in \ob J$ we define
  \[
  ((f_{\mydot})_* F^{\mydot})^j = \varprojlim (Y_\phi)_* ({f_i}_* F^i)
  \]
  where the inverse limit traverses all pairs $(i, \phi)$ where $\phi : f(i)
  \rightarrow j$ is a morphism in $J^{\opp}$.
\end{definition}

\begin{remark}
In many applications, $J$ will simply be the category $\{0 \}$  with one object and one morphism (for example, cubic hyperresolutions of schemes).  In that case one can merely think of the limit as traversing $I$.
\end{remark}

\begin{remark}
One can also define a functor $f^*$, show that it has a right adjoint and that that adjoint is $f_*$ as defined above \cite[I, 5.5]{GNPP88}.
\end{remark}

\begin{definition} \cite[I, Section 5]{GNPP88} Let $X_{\mydot}$ and $Y_{\mydot}$ be
  diagrams of topological spaces over $I$ and $J$ respectively, $\Phi : I \rightarrow
  J$ a functor, $f_{\mydot} : X_{\mydot} \rightarrow Y_{\mydot}$ a $\Phi$-morphism of
  topological spaces.  If $G^\mydot$ is a sheaf over $Y_{\mydot}$ with values in a
  complete category $\sfC$, one denotes by $f_{\mydot}^* G^\cdot$ the sheaf over
  $X_{\mydot}$ defined by
\[
(f^*_{\mydot} G^{\mydot})^i = f^*_i(G^{\Phi(i)}),
\]
for all $i \in \Ob I$.  One obtains in this way a functor
\[
f^*_{\mydot} : \uFaisc(Y_{\mydot}, \sfC) \rightarrow \uFaisc(X_{\mydot}, \sfC)
\]
\end{definition}

Given an $I$-scheme $X_{\mydot}$, one can define the category of sheaves of abelian
groups ${\sf Ab}(X_\mydot)$ on $X_{\mydot}$ and show that it has enough injectives.
Next, one can even define the derived category $D^+(X_{\mydot}, {\sf Ab}(X_\mydot))$
by localizing bounded below complexes of sheaves of abelian groups on $X_{\mydot}$ by
the quasi-isomorphisms (those that are quasi-isomorphisms on each $i \in I$).  One
can also show that $(f_{\mydot})_*$ as defined above is left exact so that it has a
right derived functor $\myR (f_{\mydot})_*$ \cite[I, 5.8-5.9]{GNPP88}.  In the case
of a cubic hyperresolution of a scheme $f : X_{\mydot} \rightarrow X$,
\[
\myR ((f_{\mydot})_* F^{\mydot}) = \myR \varprojlim (\myR {f_i}_* F^i)
\]
where the limit traverses the category $I$ of $X_{\mydot}$.

\begin{demo*}{Final Remark}%
  We end our excursion into the world of hyperresolutions here. There are many other
  things to work out, but we will leave them for the interested reader.  Many
  ``obvious'' statements need to be proved, but most are relatively straightforward
  once one gets comfortable using the appropriate language. For those and many more
  statements, including the full details of the construction of the Du~Bois complex
  and many applications, the reader is encouraged to read \cite{GNPP88}.
\end{demo*}


\def\cprime{$'$} \def\polhk#1{\setbox0=\hbox{#1}{\ooalign{\hidewidth
  \lower1.5ex\hbox{`}\hidewidth\crcr\unhbox0}}} \def\cprime{$'$}
  \def\cprime{$'$} \def\cprime{$'$} \def\cprime{$'$}
  \def\polhk#1{\setbox0=\hbox{#1}{\ooalign{\hidewidth
  \lower1.5ex\hbox{`}\hidewidth\crcr\unhbox0}}} \def\cdprime{$''$}
  \def\cprime{$'$} \def\cprime{$'$}
\providecommand{\bysame}{\leavevmode\hbox to3em{\hrulefill}\thinspace}
\providecommand{\MR}{\relax\ifhmode\unskip\space\fi MR}
\providecommand{\MRhref}[2]{%
  \href{http://www.ams.org/mathscinet-getitem?mr=#1}{#2}
}
\providecommand{\href}[2]{#2}

\end{document}